\documentclass[12pt]{amsart}
 \usepackage{amssymb, amsmath,amsfonts,amsthm,amscd,textcomp,times,mathtools }
  \ifx\pdfoutput\undefined
   \usepackage[dvips]{graphicx}
   \else
   \usepackage[pdftex]{graphicx}
   \usepackage[all,cmtip]{xy}
\usepackage[T1]{fontenc}
\usepackage{calligra}
\usepackage{multido}
\usepackage{frcursive}
\usepackage{mathrsfs}
\usepackage{tikz}
  \usepackage{float}
   \fi
  \usepackage{epstopdf}
\usepackage{cite}
 \newtheorem{theorem}{Theorem}

\newtheorem{proposition}{Proposition}
\newtheorem{lemma}{Lemma}
\newtheorem{definition}{Definition}
\newtheorem{corollary}{Corollary}
\newtheorem{remark}{Remark}
\newtheorem{example}{Example}

\numberwithin{equation}{section}
\numberwithin{theorem}{section}
\numberwithin{proposition}{section}
\numberwithin{lemma}{section}
\numberwithin{claim}{section}
\numberwithin{corollary}{section}
\numberwithin{table}{section}

\newcommand{\bull}{\ensuremath{{}\bullet{}}}
 
\newcommand{\cpn}{\ensuremath{\mathbb{P}^{N}}}
\newcommand{\slnc}{\ensuremath{SL(N+1,\mathbb{C})}}
\newcommand{\dlb}{\ensuremath{\overline{\partial}}}
\newcommand{\dl}{\ensuremath{\partial}} 
\newcommand{\ra}{\ensuremath{\longrightarrow}}
\newcommand{\om}{\ensuremath{\omega}}
\newcommand{\vp}{\ensuremath{\varphi}}
\newcommand{\vps}{\ensuremath{\varphi_{\sigma}}}

\newcommand{\mubull}{\ensuremath{\mu_{\bull}} }

\newcommand{\emu}{\ensuremath{\mathbb{E}_{\mu_{\bull}}}}
 \DeclarePairedDelimiter\abs{\lvert}{\rvert}
\makeatletter
\let\oldabs\abs
\def\abs{\@ifstar{\oldabs}{\oldabs*}}
 
\begin{document}
  \DeclareGraphicsExtensions{.pdf,.gif,.jpg}
%%%%%%%%%%%%%%%%%%%%%%%%%%%%%%%%%%%%%%%%%%%%%%%%%%%%%%%%%%%%%%%%%%%%%%%%%%%%%%%%%%%%%%%%%%%%%%%%%%%%%%%%%%%%%%%%%%
\title[Stable Pairs]{ Mahler measures, Stable Pairs, and the Global coercive estimate for the Mabuchi Functional}
\author{Sean Timothy Paul  }
\email{stpaul@wisc.edu}
\address{Mathematics Department at the University of Wisconsin, Madison}
 %%%%%%%%%%%%%%%%%%%%%%%%%%%%%%%%%%%%%%%%%%%%%%%%%%%%%%% %%%%%%%%%%%%%%%%%%%%%%%%%%%%%%%%%%%%%%%%%%%%%%%%%%
\date{May 2 , 2021} 
 \vspace{-5mm}  
\begin{abstract} We show that the Mabuchi energy of any polarized manifold $(X,L)$ is (bounded below) proper on the full space of K\"ahler metrics in the class $c_1(L)$ if and only if $(X,L)$ is asymptotically (semi)stable. It now follows from work of Xiuxiong Chen and Jinguri Cheng that $X$ admits a cscK metric in ${c}_1(L)$ iff
$(X,L)$ is asymptotically stable, provided that the group ${\tt{Aut}}(X,L)$ is finite.
%and is uniformly bounded below on the space of metrics if and only if $(X,L)$ is asymptotically semi-stable. Combining these facts with work of Xiuxiong Chen and Jinguri Cheng we conclude that $(X,L)$ admits a cscK metric in $c_1(L)$ if and only if it is stable , provided that the automorphism group $\mbox{Aut}(X,L)$ is finite.
 %We also prove the equivalence between global lower bounds for the Mabuchi energy and asymptotic semi-stability under the same hypothesis. 
\end{abstract}
\maketitle 
%\tableofcontents  
 %\newpage
%%%%%%%%%%%%%%%%%%%%%%%%%%%%%%%%%%%%%%%%%%%%%%%%%%%%%%%%%%
 %%%%%%%%%%%%%%%%%%%%%%%%%INTRO%%%%%%%%%%%%%%%%%%%%%%%%%

\section{Statement of Main Results}
 Let $G$ be a reductive algebraic group over $\mathbb{C}$. Let $\mathbb{W}$ be any finite dimensional complex representation of $G$. Fix $w\in\mathbb{W}\setminus \{0\}$ .  Define $\mathcal{O}_{w}:=G\cdot w$, the $G$ orbit of $w$ in $\mathbb{W}$. Recall that $w$ is \emph{semistable} if and only if
 \begin{align}
 0\notin \overline{\mathcal{O}}_{w}\ 
 \end{align}
 where $\overline{\mathcal{O}}_{w}$ is the Zariski closure of the orbit in $\mathbb{W}$. Next choose \emph{any} Hermitian norm $h=||\cdot ||$ on $\mathbb{W}$.  We define
 \begin{align}
 {\tt{dist}}_{h}(\overline{\mathcal{O}}_{w} , 0):=\inf \{||\sigma\cdot w|| \ |\ \sigma\in G\}\ .
 \end{align}
 Then we have the following well known proposition.
 \begin{proposition}\label{classicalsemistability}
  A point $w\in\mathbb{W}\setminus \{0\}$ is semistable if and only if there is a constant\newline $C=C(h)\geq 0 $ such that
\emph{ \begin{align}
 \log{\tt{dist}}_{h}(\overline{\mathcal{O}}_{w} , 0)\geq -C \ .
 \end{align}}
 \end{proposition}
%%%%%%%%%%%%%%%%%%%%%%%%%%%%%%%%%%%%%%%%%%%%%%%%%%%%%%%%%%%
 Let $(X,L)$ be a polarized manifold. Fix a large $k$ embedding of $X$ into $\cpn$. Let $R_X$ and $\Delta_X$ denote Cayley's  $X$-\emph{resultant} and $X$-\emph{hyperdiscriminant} respectively. Recall that these are irreducible polynomials
  in the following $G$ modules \footnote{In this paper we always take $G=\slnc$ .}
 \begin{align}
 \begin{split}
 &R_X\in  \mathbb{C}_{d(n+1)}[M_{(n+1)\times (N+1)}]^{SL(n+1,\mathbb{C})} \\
 &\Delta_X \in \mathbb{C}_{n(n+1)d-d\mu}[M_{n\times (N+1)}]^{SL(n,\mathbb{C})} \ .
 \end{split}
 \end{align}
Let $\mathcal{O}_{R\Delta}$ and $\mathcal{O}_{R}$ denote the \emph{projective} orbits 
\begin{align}
\begin{split}
&\mathcal{O}_{R\Delta}:=G\cdot [(R_X^{\deg(\Delta_X)} ,\Delta_X^{\deg(R_X)})] \subset \mathbb{P}(\mathbb{V}\oplus\mathbb{W}) \\&\mathcal{O}_{R_X}:=G\cdot [(R^{\deg(\Delta)} , 0)] \subset \mathbb{P}(\mathbb{V}\oplus\{0\}) 
\end{split}
\end{align}
where $\mathbb{V}$ and $\mathbb{W}$ are the obvious $G$ modules. Next choose \emph{any} Hermitian metric $h$ on $L$ with positive curvature. With Proposition (\ref{classicalsemistability}) in mind we make the following definition.
\begin{definition}
A polarized manifold $(X,L)$ is asymptotically semistable if and only if there is a constant  $C=C(h)\geq 0$ such that for all $k>>0$ we have
\begin{align}
 \log\tan{\tt{dist}}_{0}(\overline{\mathcal{O}}_{R\Delta} , \overline{\mathcal{O}}_{R} )\geq -Ck^{2n} \ .
\end{align}
\end{definition}
As we will explain in the sections that follow, ${\tt{dist}}_{0}$ is simply the distance between the orbit closures measured in the \emph{Mahler metric} on polynomials. The curious appearance of $\tan$ in the above definition is due to the fact that the orbits are \emph{projective}, not affine. Moreover, $R$ and $\Delta$ must be scaled to unit length in the Mahler measure. There is a similar but slightly more technical definition of asymptotic stability of $(X,L)$. This is described in detail below. For the moment we remark to the reader that any stable $X\subset \cpn$ has finite autmorphism group.

With this said we can state the main result of this article.
%%%%%%%%%%%%%%%%%%%%%%%%%%%%%%%%%%%%%%%%%%%%%%%%%%%%%%%%%%%
\begin{theorem}\label{maintheorem} Let $(X,L)$ be a polarized manifold. Let $h$ be a Hermitian metric on $L$ with positive curvature $\om_h$.   Then
\begin{itemize}
\item $(X,L)$ is asymptotically stable if and only if the Mabuchi energy is proper on $\mathcal{H}_{\om}$. \\
\item $(X,L)$ is asymptotically semistable if and only if the Mabuchi energy is bounded below on $\mathcal{H}_{\om}$. 
\end{itemize}
\end{theorem}

A variational characterization of the existence of a K\"ahler Einstein metric on a Fano manifold is provided by the following theorem of Gang Tian \cite{tian97} .\\
\ \\
\indent $(*)$ \emph{ Let $(X,\om) $ be a Fano manifold with $[\om]= c_1(X)$. Assume that $\mbox{\tt{Aut}}(X)$ is finite. Then $X$ admits a K\"ahler Einstein metric if and only if the Mabuchi energy is proper.}\\

Combining $(*)$  with our main result gives our first corollary.
 
  \begin{corollary}Let $(X, {-K_X})$ be an anti-canonically polarized manifold.  Assume that ${\tt{Aut}}(X)$ is finite. Then $(X,{-K_X})$ is asymptotically stable if and only if $X$ admits a K\"ahler Einstein metric in the class $c_1(X)$.
\end{corollary}

This provides another algebraic characterization of the existence of a K\"ahler Einstein metric on a Fano manifold with finite symmetry group.

An important development in K\"ahler geometry is the following deep result of Jinguri Cheng and Xiuxiong Chen \cite{chencheng1}, \cite{chencheng2}, \cite{chencheng3} , which generalizes Tian's properness Theorem to any K\"ahler class.\\
 
 \noindent $(**)$  \emph{Let $(X,\om)$ be a compact K\"ahler manifold. Then the Mabuchi energy is proper (modulo automorphisms of $X$, if any) on $\mathcal{H}_{\om}$ if and only if there is a metric of constant scalar curvature in the class $[\om]$.}
 \ \\

Combining $(**)$  with our main result gives our second corollary.
 \begin{corollary}Let $(X,L)$ be an arbitrary polarized manifold.  Assume that ${\tt{Aut}}(X,L)$ is finite. Then $(X,L)$ is asymptotically stable if and only if there is a constant scalar curvature metric in $c_1(L)$ .
\end{corollary}

 \subsection{Discussion}
The most important statements in the article are (\ref{distanceformula}) and (\ref{coercive}). These follow at once from Theorem \ref{precisetheorem}, the purpose of which is to \emph{identify the norms} appearing in Theorem A from \cite{paul2012}. In principle the main results of this article were available shortly after the appearance of \cite{paul2012}.  The norms that appeared in Theorem A of \cite{paul2012} are conformally equivalent to the standard $L^2$ norms on polynomials. Since the conformal factors are continuous, they are bounded by reasons of compactness. The conclusion was that the Mabuchi energy is \emph{almost} the distance between the orbits. That is, the distance in the usual Fubini Study metric induced by $L^2$ up to some (unknown) error that depended (somehow) on the degree of the embedding. Based on work in \cite{bgs1}, \cite{bgs2}, and \cite{bgs3} the author recently found a more sophisticated path to the relationship between the Mabuchi energy restricted to the Bergman metrics and the resultant and hyperdiscriminant of the subvariety which revealed that the error was in fact the {difference} between the $L^2$ norm and \emph{another} well known norm, namely the $L^0$ norm\footnote{ That the $L^p$ norms only give norms for $p\in[1,\infty]$ does not matter in this article.}, i.e. the Mahler measure (see Theorem (\ref{precisetheorem})) .
The boundedness of the error, initially attributed to compactness, is just an expression of the fact that these norms are comparable.  The outcome is that the norm on the space of polynomials which connects the Mabuchi energy to stability of the pair $(R,\Delta)$ is \emph{exactly} given by the Mahler measure.  Now asymptotic stability and global bounds on K-energy maps follow almost at once from Tian's Thesis \cite{tianberg} .

   The strategy of restricting to the Bergman metrics is due to Tian and Yau. Tian explained it to the author many years ago. Despite  Tian's many works on the subject, as well as the articles \cite{donaldson2001} and \cite{donaldson2005}, this strategy was never really developed.  Instead, the approach of Tian in \cite{tianCPAM} as well as Chen-Donaldson-Sun in \cite{cds1},\cite{cds2},\cite{cds3} is to reduce an infinite dimensional estimate to a finite dimensional one. Whereas the approach of this author is to obtain the infinite dimensional estimate from a sequence of finite dimensional estimates. \begin{center}\emph{The finite dimensional estimates are equivalent to the stability of the variety with respect to the given embedding.}\end{center}
 As we have mentioned, the precise definition of the asymptotic stability of a polarized manifold appears below. The relevant ideas are contained in definitions (\ref{semistability}),(\ref{stability}),(\ref{asymptoticsemistability}), and (\ref{asymptoticstability}). The reader should compare the author's definition of (semi)stability with the many variations of K-Stability that appear in the literature. First, we consider orbits under all of $G$, not just one parameter subgroups of $G$. Second, from the author's point of view, stability is not necessarily concerned with a variety in a projective space. Stability is a property that a pair of (non-zero) vectors in a pair of finite dimensional complex representations of an algebraic group may, or may not, possess. As the reader shall see, the stability of a projective variety is a special case of this situation. As we have already mentioned, {test configurations} do not play a direct part in our definition of stability, they are rather considered as a \emph{means to check} stability.  This is exactly how one parameter subgroups are used in Hilbert and Mumford's Geometric Invariant Theory. The author is optimistic about the eventual conversion of our stability condition into a (hopefully tractable) combinatorial condition that can be checked in concrete examples. This is due to the fact that a (semi)stable pair is a straightforward extension of Mumford's stability (see the table at the end of section 2) . In particular the vast array of tools concerning actions of reductive groups on finite dimensional representations can be used. On the other hand, the author does not expect that checking the stability of a given pair will ever be made entirely trivial:  indeed, even after so many years, and so much creative effort, checking the stability of Chow and Hilbert points in dimensions $\geq 3$ still seems to be out of reach. 

This paper is organized as follows. In section 2 we give an account of the (semi)stability of pairs of rational representations of a reductive algebraic group and provide several elementary examples of such pairs, we also show that the automorphism group of a stable pair is finite. In section 3 we apply the ideas of section 2 to (hyper)discriminants and resultants of a projective variety. In section 4 we discuss the equivalence among the $L^p$ metrics on orbit closures in polynomial representations. Of special importance is the case when $p=0$. The special metric induced by this measure allows us to define asymptotic (semi)stability for any polarized manifold $(X,L)$.  
In section 5 we write down the conformal factor that appears in Theorem A of \cite{paul2012}. This allows us to show that the Mabuchi energy \emph{is} the distance between the orbits where the distance is computed the the Mahler metric. This is enough to establish the equivalence between the global coercive estimate for the Mabuchi energy and the asymptotic stability of the polarized manifold.    
%%%%%%%%%%%%%%%%%%%%%%%%%%%%%%%%%%%%%%%%%%%%%%%%%%%%%% 
 \section{Semistability of Pairs}
 Let $G$ denote any of the classical linear reductive algebraic groups over $\mathbb{C}$.  For example, $G$ can be taken to be any one of the following
\begin{align}
GL(N,\mathbb{C}), SL(N,\mathbb{C}) \ , \ SO(N,\mathbb{C}) \ ,\  O(N,\mathbb{C})\ , \ Sp(2N,\mathbb{C}) \ .
\end{align}
Primarily we will be interested in the case when $G$ is the special linear group.
      For any vector space $\mathbb{V}$ and any $v\in \mathbb{V}\setminus\{0\}$ we  let $[v]\in\mathbb{P}(\mathbb{V})$ denote the line through $v$. If $\mathbb{V}$ is a $G$ module then we can consider the \emph{\textbf{projective}} orbit :     
\begin{align}
  \mathcal{O}_{v}:=G\cdot [v]  \subset \mathbb{P}(\mathbb{V})\ .
\end{align}
We let $ \overline{\mathcal{O}}_{v}$ denote the Zariski closure of this orbit.  

We consider pairs $(\mathbb{E}; e)$ such that $\mathbb{E}$ is a finite dimensional complex $G$-module and the linear span of the orbit $G\cdot e$ coincides with $\mathbb{E}$ .  
    
\begin{definition} (see \cite{smirnov2005})   
A pair $(\mathbb{U}; u)$ \textbf{{dominates}} $(\mathbb{W}; w)$, in which case we write $(\mathbb{U}; u)\succsim (\mathbb{W}; w)$,  
 if and only if there exists $\pi\in Hom(\mathbb{U},\mathbb{W})^G$ such that
$ \pi(u)=w$ and the induced  rational map  
$\pi:\mathbb{P}(\mathbb{U}) \dashrightarrow  \mathbb{P}(\mathbb{W})$
restricts to a regular finite map  
$\pi:\overline{\mathcal{O}}_{u}\ra \overline{\mathcal{O}}_{w} \ $  between the Zariski closures of the orbits.
  \end{definition} 
 {My approach to the {Stability Conjectures} is based on this definition. In \cite{smirnov2005}, the motivation behind making such a definition seems to the problem of decomposing the {symmetric} power of an irreducible representation of $GL(n,\mathbb{C})$. It is mysterious that the \emph{same} definition appears\footnote{The author was led to the same definition independently. See ``semistable pair'' below.} when one seeks to bound (from below) the Mabuchi energy restricted to the space of Bergman metrics.}
 
Observe that the restriction of the map $\pi$  to $\overline{\mathcal{O}}_{u}$ is regular if and only if the following holds
  \begin{align}\label{meetskernel}
 \qquad  \overline{\mathcal{O}}_{u}\cap \mathbb{P}({\tt{ker}}(\pi))=\emptyset \ .
\end{align}
As the reader can easily check, whenever $(\mathbb{U}; u)\succsim (\mathbb{W}; w)$ it follows that
\begin{align} 
 & \pi(\mathbb{U})=\mathbb{W} \ \mbox{and} \ \mathbb{U}={\tt{ker}}(\pi)\oplus \mathbb{W} \ \mbox{ ($G$-module splitting) } \ .
\end{align}
Therefore we may identify $\pi$ with projection onto $\mathbb{W}$ and $u$ decomposes as follows
\begin{align}
v=(u_{\pi},w) \ , \ {\tt{ker}}(\pi)\ni u_{\pi}\neq 0 \ .
\end{align}
Again the reader can easily check that (\ref{meetskernel}) is equivalent to
\begin{align}
 \qquad \overline{\mathcal{O}}_{(u_{\pi},w)}\cap\overline{\mathcal{O}}_{u_{\pi}}=\emptyset \quad \mbox{( Zariski closure in  $\mathbb{P}({\tt{ker}}(\pi)\oplus\mathbb{W}$ ) )} \ .
\end{align}
   
We summarize this discussion in the following way. Given $\mathbb{V}$ and $\mathbb{W}$ two $G$ representations with (nonzero) points $v$ and $w$ respectively, we consider, as before, the projective orbits\footnote{ We do not assume anything about the linear spans of the orbits. }     
\begin{align}
\mathcal{O}_{(v,w)}:=G\cdot [(v,w)]  \subset \mathbb{P}(\mathbb{V}\oplus\mathbb{W}) \ , \ \mathcal{O}_{v}:=G\cdot [(v,0)]  \subset \mathbb{P}(\mathbb{V}\oplus\{0\})\ .
\end{align}
  Now we can give the definition of a \emph{semistable pair}. This definition seems the most appropriate for the Stability Conjectures as it gives precise characterization of the infimum of the Mabuchi energy restricted to the space of Bergman metrics.  
   \begin{definition}\label{semistable} 
 {The pair $(v,w)$ is \textbf{ {semistable}} if and only if}  $ \overline{\mathcal{O}}_{(v,w)}\cap\overline{\mathcal{O}}_{v}=\emptyset $ .
  \end{definition} 
The relationship of this with Mumford's Geometric Invariant Theory is brought out in the following example.
 
 \begin{example} 
\emph{ Let $\mathbb{V}\cong \mathbb{C}$ be the trivial one dimensional representation and let $v=1$ . Suppose $\mathbb{W}$ is any representation of $G$ and let $w\in \mathbb{W}\setminus\{0\}$ . Then 
 $([1],[w])$ is a semistable pair if and only if $0\notin \overline{G\cdot w}\subset \mathbb{W}$ .}
 \end{example}  
\begin{example} \emph{Let $\mathbb{V}_e$ and  $\mathbb{V}_d$ be irreducible $SL(2,\mathbb{C})$ modules with highest weights $e\in\mathbb{N}$ and $d\in\mathbb{N}$ respectively. These are well known to be spaces of homogeneous polynomials in two variables. Let $f$ and $g$ be two such polynomials in $\mathbb{V}_e\setminus\{0\}$ and $\mathbb{W}_d\setminus\{0\}$
respectively. If the pair $(f,g)$ is semistable  then we must have that
\begin{align}
e\leq d \ \mbox{and for all $p\in \mathbb{P}^1$ we have }\  {\tt{ord}}_p(g)-{\tt{ord}}_p(f)\leq \frac{d-e}{2} \ .
\end{align}}
\emph{In particular when $e=d-1$ \emph{there are no semistable pairs} .} 
\end{example}

 Let $\mathbb{E}$ be a finite dimensional {reducible} representation of $G$. Let $u\in \mathbb{E}\setminus\{0\}$. Let $\mathcal{O}\subset \mathbb{P}(\mathbb{E})$ denote the projective orbit $G\cdot [u]$. We assume that the linear span of $\mathcal{O}$ coincides with $\mathbb{P}(\mathbb{E})$ . Fix a Borel subgroup $B\leq G$ and a maximal algebraic torus $T\leq B$. Let $\Lambda^+$ denote the dominant integral weights relative to $B$. It is well known that  $\overline{\mathcal{O}}$ is a union of orbits at least one of which is {closed} and each closed orbit corresponds to an irreducible $G$ submodule $\emu$ of $\mathbb{E}$. We assume that $\overline{\mathcal{O}}$ consists of \emph{finitely many} orbits. Let $\Lambda^+(\mathcal{O})$ denote the dominant weights corresponding to the closed orbits in $\overline{\mathcal{O}}$. Then we have the decomposition
\begin{align*}
\overline{\mathcal{O}}=\mathcal{O} \cup \mathcal{O}_{1}\cup \dots \cup  \mathcal{O}_{k} \bigcup _{\mubull\in \Lambda^+(\mathcal{O})}G\cdot[w_{\mubull}] \ ,
\end{align*}
where $w_{\mubull}$ is the corresponding highest weight vector . Now we decompose $\mathbb{E}$ according to the orbit $\mathcal{O}$
\begin{align*}
\mathbb{E}=  \mathbb{V} \oplus \bigoplus _{\mubull\in \Lambda^+(\mathcal{O})}\emu   \ .
\end{align*}
  We assume that $\mathbb{V}\neq 0$ . Let $\pi_{\mathcal{O}}$ and $\pi_{\mathbb{V}}$ denote the projections onto $\bigoplus _{\mubull\in \Lambda^+(\mathcal{O})}\emu$ and $\mathbb{V}$ respectively.  Then we may decompose $u$ as follows
  \begin{align*}
  u= (v,w):= (\pi_{\mathbb{V}}(u),\pi_{\mathcal{O}}(u)) \ .
  \end{align*}
Then $(v,w)$ is semistable if and only if  for every $1\leq i\leq k$ there exists a $\mubull\in \Lambda^+(\mathcal{O})$  such that $\pi_{\mubull}(x_i)\neq 0$  where $\mathcal{O}_i=G\cdot [x_i]$ and 
$\pi_{\mubull}$ is the projection onto $\emu$.

The simplest situation is the case when $\overline{\mathcal{O}}$ consists of {two} orbits (one of which is closed)
\begin{align*}
\overline{\mathcal{O}}= \mathcal{O}\cup G\cdot [w_{\mubull}] \ .
\end{align*}
In this case it is automatic that the pair $(\pi_{\mathbb{V}}(u),\pi_{\mathcal{O}}(u))$ is semistable. Therefore the class of two orbit varieties ( or, more generally, \emph{quasi-closed orbits} ) provides many examples of semistable pairs. Such varieties have been completely classified by Stephanie Cupit-Foutou (see \cite{cupit2003}) and Alexander Smirnov (see \cite{smirnov2004}).
%%%%%%%%%%%%%%%%%%%%%%%%%%%%%%%%%%%%%%%%%%%%%%%%%%%%%%%%%%%%%%%%%%%%%%%%
 
%%%%%%%%%%%%%%%%%%%%%%%%%%%%%%%%%%%%%%%%%%%%%%%%%%%%%%
\begin{example}\label{blowup}\emph{
 Let $\psi:\mathbb{P}^2\times\mathbb{P}^2\dashrightarrow \mathbb{P}(\wedge^2\mathbb{C}^3)$ be the rational map $\psi([v],[w]):=[v\wedge w]$.
  The graph of $\psi$ is
\begin{align*}
\Gamma_{\psi}:=\{([v],[w],[v\wedge w])\ |\ [v]\neq [w]\}\subset \mathbb{P}^2\times\mathbb{P}^2\times\mathbb{P}(\wedge^2\mathbb{C}^3) \ .
\end{align*}
Recall that the blow up of $\mathbb{P}^2\times\mathbb{P}^2$ along the diagonal $\Delta$  is the Zariski closure of $\Gamma_{\psi}$ inside $\mathbb{P}^2\times\mathbb{P}^2\times\mathbb{P}(\wedge^2\mathbb{C}^3)$ . We will denote the blow up by  $ B_{\Delta}(\mathbb{P}^2\times\mathbb{P}^2)$ and let $E\cong \mathbb{P}(T^{1,0}_{\mathbb{P}^2})$ denote the exceptional divisor. The situation can be pictured as follows
\begin{align}
\begin{split}
\xymatrix{
 \Gamma_{\psi}\subset B_{\Delta}(\mathbb{P}^2\times\mathbb{P}^2)\ar@{^{(}->}[r]^-{\iota} \ar[d]^{p_{12}} & \mathbb{P}^2\times\mathbb{P}^2\times\mathbb{P}(\wedge^2\mathbb{C}^3)   \ar[r]^-{S}\ar[d]^{p_3}& \mathbb{P}(\mathbb{E}_{310}\oplus \mathbb{C}^3\oplus S^2(\wedge^2\mathbb{C}^3)\oplus  \mathbb{C}^3)   \\
\mathbb{P}^2\times\mathbb{P}^2\ar@{-->}^{\psi}[r]&\mathbb{P}(\wedge^2\mathbb{C}^3)  } \ .
\end{split}
\end{align}
Then we claim that $B= B_{\Delta}(\mathbb{P}^2\times\mathbb{P}^2)$ is a two-orbit $G=SL(3,\mathbb{C})$ variety (for the natural $G$ action) with  orbit decomposition
\begin{align}
B=(B\setminus E) \cup E \ .
\end{align}
Where $(B\setminus E)$ is necessarily the open orbit. There is a $G$ equivariant identification
\begin{align*}
B\setminus E\cong \mathbb{P}^2\times\mathbb{P}^2\setminus \Delta \ .
\end{align*}
Since $G$ acts transitively on planes in $\mathbb{C}^3$ we easily get that $\mathbb{P}^2\times\mathbb{P}^2\setminus \Delta$ is an orbit:
\begin{align*}
G\cdot ([e_1],[e_2]) = \mathbb{P}^2\times\mathbb{P}^2\setminus \Delta \ .
\end{align*}
To see that $E$ is a homogeneous $G$ variety we can proceed as follows. We have the decomposition into  irreducible summands
\begin{align*}
\mathbb{C}^3\otimes \mathbb{C}^3\otimes \wedge^2\mathbb{C}^3\cong \mathbb{E}_{310}\oplus \mathbb{C}^3\oplus S^2(\wedge^2\mathbb{C}^3)\oplus  \mathbb{C}^3 \ .
\end{align*}
The summand $\mathbb{E}_{310}$ appears as follows
\begin{align*}
0\ra \mathbb{E}_{310}\cong\mbox{Ker}(\pi)\ra S^2(\mathbb{C}^3)\otimes \wedge^2(\mathbb{C}^3)\xrightarrow{\pi} \mathbb{C}^3\ra 0 \ ,
\end{align*}
where the map $\pi$ is defined by 
\begin{align*}
\pi(v\cdot w\otimes \alpha)=\alpha(v)w+\alpha(w)v \ .
\end{align*}
Note that
\begin{align*}
e_1^2\otimes (e_1\wedge e_2) \in \mbox{Ker}(\pi) \ .
\end{align*}
Since $e_1^2\otimes (e_1\wedge e_2)$ is a highest weight $(310)$ vector we see that $\mathbb{E}_{310}$ is  a summand of  $\mbox{Ker}(\pi)$.  Since these spaces have the same dimension (which is 15 by the Weyl dimension formula) they coincide. Next we observe that
\begin{align*}
([e_1+te_2], [e_1], [e_1\wedge e_2])\in \Gamma_{\psi} \quad \mbox{for all $t\in\mathbb{C}^*$}\ . 
\end{align*}
As $t\ra 0$ we have
\begin{align*}
([e_1+te_2], [e_1], [e_1\wedge e_2])\ra([e_1],[e_1],[e_1\wedge e_2]) \in E \ .
\end{align*}
Let $S:\mathbb{P}^2\times\mathbb{P}^2\times\mathbb{P}(\wedge^2\mathbb{C}^3)\ra \mathbb{P}(\mathbb{E}_{310}\oplus \mathbb{C}^3\oplus S^2(\wedge^2\mathbb{C}^3)\oplus  \mathbb{C}^3)$ denote the 
Segre map. Then we have that
\begin{align*}
S([e_1],[e_1],[e_1\wedge e_2])=[e_1^2\otimes (e_1\wedge e_2)] \ .
\end{align*}
Therefore 
\begin{align*}
S(E)= G\cdot [e_1^2\otimes (e_1\wedge e_2)]  \ .
\end{align*}
Since $S$ is an embedding $E$ is a closed orbit with stabilizer
\begin{align*}
\begin{pmatrix}
*&*&*\\
0&*&*\\
0&0&*
\end{pmatrix}
\end{align*}
therefore we identify $E$ with ${F}(1,2,\mathbb{C}^3)$ the space of complete flags in $\mathbb{C}^3$. 
The projection
\begin{align*}
{F}(1,2,\mathbb{C}^3)\xrightarrow{p_1}\mathbb{P}^2 \ 
\end{align*}
exhibits ${F}(1,2,\mathbb{C}^3)$ as a projective bundle with fiber
\begin{align*}
p_1^{-1}([v])=\mathbb{P}(\mathbb{C}^3/\mathbb{C}v)  \ .
\end{align*}
Therefore if $\mathcal{Q}$ denotes the quotient bundle over $\mathbb{P}^2$ then we have the $G$ equivariant identifications 
\begin{align*}
{F}(1,2,\mathbb{C}^3)\cong \mathbb{P}(\mathcal{Q})\cong \mathbb{P}(\mathcal{O}(1)\otimes\mathcal{Q}) = \mathbb{P}(T^{1,0}_{\mathbb{P}^2})
\end{align*}
as expected. $S$ maps the point $[(v\otimes w\otimes (v\wedge w)]$ in $X\setminus E$ to
\begin{align*}
v\cdot w\otimes (v\wedge w)+(v\wedge w)^2 \in \mathbb{E}_{310}\oplus S^2(\wedge^2\mathbb{C}^3) \cong \mathbb{E}_{310}\oplus \mathbb{E}_{220}\ .
\end{align*}
We conclude that the pair
\begin{align*}
(e,f):=\big( (e_1\wedge e_2)^2 \ , \ e_1\cdot e_2\otimes (e_1\wedge e_2)\big)\in \mathbb{E}_{220}\oplus \mathbb{E}_{310} 
\end{align*}
is semistable.
 }
\end{example}
$\Box$
\begin{remark}
\emph{The semistability of the pair $(v,w)$ depends only on $([v],[w])$. The reader should also observe that the definition is \emph{not} symmetric in $v$ and $w$. In virtually all examples where the pair $(v,w)$ is semistable $(w, v)$ is not semistable.}
\end{remark}
\subsection{Numerical Semistability}
If the pair $(v,w)$ is semistable then obviously we have 
\begin{align}\label{numerical}
\overline{T\cdot[(v,w)]}\cap \overline{T\cdot[(v,0)]}=\emptyset 
\end{align}
for all algebraic tori $T$ of $G$. We may as well assume that $T$ is \emph{maximal}. In this section we relate semistability to lattice polytopes.  To begin we let $M_{\mathbb{Z}} $ be the {character lattice} of $T$
\begin{align}
M_{\mathbb{Z}}:= \mbox{Hom}_{\mathbb{Z}}(T,\mathbb{C}^*) \ . 
\end{align}

As usual, the dual lattice is denoted by $N_{\mathbb{Z}}$. It is well known that $ u\in N_{\mathbb{Z}}$ corresponds to an algebraic one parameter subgroup $\lambda$ of $T$. These are algebraic homomorphisms $$\lambda:\mathbb{C}^*\ra T\ .$$ The correspondence is given by 
\begin{align}
(\cdot \ , \ \cdot) :N_{\mathbb{Z}}\times M_{\mathbb{Z}}\ra \mathbb{Z} \ , \ m(\lambda(\alpha))=\alpha^{(u , m)} \quad m\in M_{\mathbb{Z}}\ .
\end{align}
 
 We introduce associated real vector spaces by extending scalars 
\begin{align} 
 \begin{split}
 M_{\mathbb{R}}:= M_{\mathbb{Z}}\otimes_{\mathbb{Z}}\mathbb{R} \qquad N_{\mathbb{R}}:= N_{\mathbb{Z}}\otimes_{\mathbb{Z}}\mathbb{R}\ .
\end{split}
\end{align}
Then the one parameter subgroups $\lambda$ of $T$ may be viewed as integral linear functionals $$l_{\lambda}:M_{\mathbb{R}}\ra \mathbb{R}\ . $$
Any rational representation $\mathbb{E}$ decomposes under the action of $T$ into  {weight spaces}
\begin{align}
\begin{split}
\mathbb{E}=\bigoplus_{a\in {\mathscr{A} }}\mathbb{E}_{a}  \qquad \mathbb{E}_{a}:=\{e\in \mathbb{E}\ |\ t\cdot e=a(t) e \ , \ t\in T\}
\end{split}
\end{align}
$\mathscr{A} $ denotes the $T$-{support} of $\mathbb{E}$ 
\begin{align}
\mathscr{A}:= \{a \in M_{\mathbb{Z}}\ | \ \mathbb{E}_{a}\neq 0\} \ .
\end{align}
Observe that $\mathscr{A}$ is a finite subset of $M_{\mathbb{Z}}$.

\indent Given $e\in \mathbb{E}\setminus \{0\}$  the projection of $e$ into $\mathbb{E}_{a}$ is denoted by $e_{a}$. The support of any (nonzero) vector $e$ is then defined by
\begin{align}
\mathscr{A}(e):= \{a\in \mathscr{A}\ | \ e_{a}\neq 0\} \ .
\end{align}
\begin{definition}  \emph{ Let $T$ be any maximal torus in $G$. Let $e\in \mathbb{E}\setminus\{0\}$ . The \textbf{\emph{weight polytope}} of $e$ is the compact convex lattice polytope $\mathcal{N}(e)$ given by}
\begin{align}\label{wtpolytope}
\mathcal{N}(e):=\mbox{ {\tt{conv}}}\  \mathscr{A}(e) \subset M_{\mathbb{R}} \ 
\end{align}
\emph{where ${\tt{conv}}\mathscr{A}(e)$ denotes the convex hull of the finite set $\mathscr{A}(e)$ }.
\end{definition}
 
\begin{definition}\label{weight}
\emph{Let $\mathbb{E}$ be a rational representation of $G$. Let $\lambda$ be any 1psg of $T$  . The \textbf{\emph{weight}}  $w_{\lambda}(e)$  of $\lambda$ on $e\in \mathbb{E}\setminus\{0\}$ is the integer}
\begin{align}
w_{\lambda}(e):= \mbox{\emph{min}}_{   x\in \mathcal{N}(e) }\ l_{\lambda}(x)= \mbox{\emph{min}} \{ (a,\lambda)\ | \ a \in \mathscr{A}(e)\}\ .
\end{align}
 \noindent\emph{Alternatively, $w_{\lambda}(e)$ is the unique integer such that}
\begin{align}
\lim_{|t|\rightarrow 0}t^{-w_{\lambda}(e)}\lambda(t)e \  \mbox{ {exists in $\mathbb{E}$ and is \textbf{not} zero}}.
\end{align}
\end{definition} 
Next, given $d\in\mathbb{N}$ and $a\in\mathscr{A}$ recall that the $T$ {semi-invariants} $P\in \mathbb{C}_d{[\mathbb{E}]_a}^T$ of degree $d$
are characterized by
$$ P(\tau\cdot e)=a(\tau)P(e)\quad \mbox{for all $\tau\in T$}\ .$$

 \begin{proposition} 
Let $T$ be a maximal algebraic torus of $G$, and let $\mathbb{V}$ and $\mathbb{W}$ two finite dimensional rational $G$-modules. Then the following are equivalent: 
\begin{align} 
\begin{split}
&1) \quad \overline{T\cdot[(v,w)]}\cap \overline{T\cdot[(v,0)]}=\emptyset \\
\vspace{4mm}
&2) \quad \mathcal{N}(v)\subset \mathcal{N}(w) \\ 
 \vspace{4mm}
&3)\quad  w_{\lambda}(w)\leq w_{\lambda}(v) \ \mbox{for all 1psg's}\ \lambda :\mathbb{C}^*\ra T \\
\vspace{4mm}
  &4) \quad \mbox{For every $\chi\in\mathscr{A}(v)$}\ \mbox{there is an}\   f\in  \mathbb{C}_d[\ \mathbb{V}\oplus\mathbb{W}\ ]^T_{d\chi} \\
  & \qquad\mbox{such that}\ f((v,w))\neq 0\ \mbox{and}\ f|_{\mathbb{V}}\equiv 0\ .
  \end{split}
    \end{align}
\end{proposition}
 \begin{proof} The equivalence of 1) and 2) follows by a simple modification of the argument in \cite{birkes71}, the equivalence of 1) and 4) follows from the Nullstellensatz, the remaining equivalences are left to the reader .\end{proof}
 
There should be an analogue of the Hilbert Mumford numerical criterion in our situation.\\
\textbf{Question.}  
\emph{ In addition to the requirement that $\mathcal{N}(v)\subset \mathcal{N}(w) $ for all maximal algebraic tori $T\leq G$ are further combinatorial conditions required to insure that the pair $(v,w)$ is semistable ?} \\

%%%%%%%%%%%%%%%%%%%%%%%%%%%%%%%%%%%%%%%%%%%%% 
In order to define a \emph{\textbf{strictly stable}} (henceforth stable) pair we need a large (but fixed) integer $m$ and the auxiliary left regular representation of $G$
\begin{align}
G\times\mathcal{GL}(N+1,\mathbb{C}) \ \ni \ (\sigma , A)\ra \sigma\cdot A \ .
\end{align}
  Recall that $\mathcal{GL}(N+1,\mathbb{C})$ is the vector space of square matrices of size $N+1$. The action is matrix multiplication. 
The {{standard $N$-simplex}}, denoted by $Q_N$, is defined to be the weight polytope of the identity operator 
\begin{align}
{I}\in \mathcal{GL}(N+1,\mathbb{C})\ .
\end{align}
 $Q_N$  is full-dimensional and contains the origin in its strict interior
\begin{align}
0\in Q_N:=\mathcal{N}( {I})\subset M_{\mathbb{R}} \ .
\end{align}  

Let $\mathbb{V}$ be a $G$ module. We define the \emph{\textbf{degree}} of $\mathbb{V}$  as follows
\begin{align}
\deg(\mathbb{V}):=\min\Big\{k\in \mathbb{Z}_{>0} \ |\ \mathcal{N}(v)\subseteq kQ_N\ \mbox{for all $0\neq v\in \mathbb{V}$} \ \Big\} \ .
\end{align}

For example, if $G=\slnc$ and $\mathbb{V}={\tt{Sym}}^d(\mathbb{C}^{N+1})^{\vee}$ then the degree of $\mathbb{V}$ is $d$.

Let $v\in \mathbb{V}$, $w\in \mathbb{W}$, and $m\in \mathbb{N}$ we define
\begin{align}
\begin{split}
& v^m:=v^{\otimes m}\in \mathbb{V}^{\otimes m} \ , \ w^{m+1}:=w^{\otimes (m+1)}\in \mathbb{W}^{\otimes (m+1)}\\
&  {I}^q:={I}^{\otimes q}\in \mathcal{GL}(N+1,\mathbb{C})^{\otimes q} \ .
\end{split}
\end{align}
Finally we can give the definition of a stable pair. 
\begin{definition}\label{stable} 
\emph{The pair $(v,w)$ is \emph{\textbf{stable}} if and only if there is a positive integer $m$ such that $({I}^q\otimes v^m \ , \ w^{m+1})$ is semistable where $q$ denotes the degree of $\mathbb{V}$. }
\end{definition}

We define the \emph{automorphism group} of the pair $(v,w)$ as
\begin{align}
{\tt{Aut}}(v,w):= G_{[v]}\cap G_{[w]} \ .
\end{align}
We have developed enough of the theory of (semi)stable pairs in this section to state the following proposition.
\begin{proposition}\label{finiteautomorphism}
The automorphism group of a stable pair is finite.
\end{proposition}
\begin{proof}
Stability of $(v,w)$ is equivalent to the inequality
\begin{align}
m(\log||\sigma\cdot w||-\log||\sigma\cdot v||) \geq  \deg(\mathbb{V})\log||\sigma||-\log||\sigma\cdot v|| \ 
\end{align}
for all $\sigma\in G$ where $m$ is a positive integer.  Decompose ${\tt{Aut}}(v,w)$  into it's reductive ($S$)  and unipotent ($U$) parts
\begin{align}
{\tt{Aut}}(v,w)=S\cdot U \ .
\end{align}
Since $U$ has no non-trivial characters stability implies that there is a constant $C$ such that
\begin{align}
C\geq \log||u|| \quad \mbox{for all $u\in U$ .}
\end{align}
Since a (euclidean) bounded affine algebraic variety is a finite collection of points, we see that $U$ must be finite. Since stability implies semistability  the weights of any $\lambda:\mathbb{C}^{*}\ra {\tt{Aut}}(v,w)$ must coincide. Precisely
\begin{align}
w_{\lambda}(w)=w_{\lambda}(v) \quad \mbox{for all 1psg's $\lambda$ of $S$} \ .
\end{align}
Once more stability shows that for all such $\lambda$ we have
\begin{align}
\deg(\mathbb{V})w_{\lambda}(I)-w_{\lambda}(v)\geq 0 \ .
\end{align}
On the other hand, by definition of the degree of a representation we have
\begin{align}
\mathcal{N}(v)\subset \deg(\mathbb{V})\mathcal{N}(I) 
\end{align}
which implies equality (remember that bigger polytopes have \emph{smaller} weights )
\begin{align}
\deg(\mathbb{V})w_{\lambda}(I)=w_{\lambda}(v)\ .
\end{align}
Since $\lambda$ lies in ${\tt{Aut}}(v,w)$ we have
\begin{align}
w_{\lambda}(v)=-w_{\lambda^{-1}}(v) \ .
\end{align}
Therefore we see that for all $\lambda$ in $S$ we have
\begin{align}\label{weightequality}
w_{\lambda}(I)=-w_{\lambda^{-1}}(I) \ .
\end{align}
Observe that we may diagonalize $\lambda$ 
\begin{align}
\lambda=(a_0\geq a_1\geq \dots \geq a_N) \quad a_i\in \mathbb{Z} \quad \sum_{0\leq i\leq N}a_i=0 \ .
\end{align}
(\ref{weightequality}) implies that $a_0=a_N$ and therefore that $\lambda$ is trivial. This completes the proof.
\end{proof}
%%%%%%%%%%%%%%%%%%%%%%%%%%%%%%%%%%%%%%%%%%%%%%%%%%%%%%%%%%%%%% 
Choose Hermitian inner products on $\mathbb{V}$ and $\mathbb{W}$ . If we give $\mathbb{V}\oplus\mathbb{W}$ the orthogonal sum metric then we may define the usual Fubini Study Riemannian metric  $g_{FS}$ on $\mathbb{P}( \mathbb{V}\oplus\mathbb{W} )$. Choose any  $\sigma , \tau \in G $. The well known formula for the distance between two points in the Fubini Study metric gives the inequality
\begin{align}
\cos{ \tt{dist}}_{g_{FS}}(\sigma\cdot[(v,w)],\tau\cdot[(v,0)])\leq \cos{ \tt{dist}}_{g_{FS}}(\sigma\cdot[(v,w)],\sigma\cdot[(v,0)])\ .
\end{align} 
In particular we see that
\begin{align} 
{ \tt{dist}}_{g_{FS}}(\sigma\cdot[(v,w)],\sigma\cdot[(v,0)])\leq { \tt{dist}}_{g_{FS}}(\sigma\cdot[(v,w)],\tau\cdot[(v,0)])\ .
\end{align}
The fact that the group elements are the same on the left hand side of this inequality implies the following.
\begin{proposition}
We have the identity
\begin{align}  
\inf_{\sigma\in G}{ \tt{dist}}_{g_{FS}}(\sigma\cdot[(v,w)],\sigma\cdot[(v,0)])= { \tt{dist}}_{g_{FS}}(\overline{\mathcal{O}}_{(v,w)} ,\overline{\mathcal{O}}_{v})\ .
\end{align}
\end{proposition}
 Another direct application of the distance formula gives
  \begin{align}
 \log \tan^2{ \tt{dist}}_{g_{FS}}(\sigma\cdot[(v,w)],\sigma\cdot[(v,0)])=\log||\sigma\cdot w||^2-\log||\sigma\cdot v||^2 \ .
  \end{align}
  
  Therefore we have the following proposition.
  \begin{proposition}\label{logdistance}
  We have the identity
  \begin{align}  
   \inf_{\sigma\in G}\left( \log||\sigma\cdot w||^2-\log||\sigma\cdot v||^2\right)=\log\tan^2 { \tt{dist}}_{g_{FS}}(\overline{\mathcal{O}}_{(v,w)} ,\overline{\mathcal{O}}_{v}) \ .
 \end{align}
  \end{proposition}
 %{\remark{This fact, while completely elementary, is what allows our proof of the Yau-Tian-Donaldson conjecture to go through. }}
   
We end this section with a direct comparison of Mumford's stability and the author's stability of pairs. Observe that the left hand column of Table \ref{mumford-me} below arises from the right when we take $\mathbb{V}\cong\mathbb{C}$ (the trivial one dimensional representation) and $v=1$. Recall that $q$ denotes the degree of $\mathbb{V}$ .
  
%%%%%%%%%%%%%%%%%%%%%%%%%%%%%%%%%%%%%%%%%%%%%%%%%%%%%%%%%%%%%%%%%%%%%
\begin{table}[h]\label{mumford-me}\caption{\sc{Hilbert Mumford Semistability vs. Semistable Pairs}} \begin{tabular}{l|l}
 \hline \rule{0pt}{3ex}
 For all $T\leq G$ $\exists\ d\in \mathbb{Z}_{>0}$ and & For all $T\leq G$ and $\chi\in\mathscr{A}(v)$\\
 $f\in \mathbb{C}_{\leq d}[\ \mathbb{W}\ ]^T$ such that   & $\exists\ d\in \mathbb{Z}_{>0}$ and $f\in  \mathbb{C}_d[\ \mathbb{V}\oplus\mathbb{W}\ ]^T_{d\chi}$  \\
   $f(w)\neq 0$ and $f(0)=0$ & such that $f((v,w))\neq 0$ and $f|_{ \mathbb{V}}\equiv 0$ \\ 
  \hline \rule{0pt}{3ex}
 $0\notin\overline{G\cdot w}$ & \ $\overline{\mathcal{O}}_{(v,w)}\cap\overline{\mathcal{O}}_{v}=\emptyset$  \rule[-1.2ex]{0pt}{0pt}  \\ 
 \hline  \rule{0pt}{3ex}
 $w_{\lambda}(w)\leq 0$   &\ $w_{\lambda}(w)-w_{\lambda}(v)\leq 0$   \\
 for all 1psg's $\lambda$ of $G$ & \ for all 1psg's $\lambda$ of $G$\\ 
 \hline \rule{0pt}{3ex}
 $0\in \mathcal{N}(w)$ all $T\leq G$ &\ $\mathcal{N}(v)\subset \mathcal{N}(w)$ {all $T\leq G$}  \rule[-1.2ex]{0pt}{0pt}\\ 
 \hline  \rule{0pt}{3ex}
 $\exists$ $C\geq 0$ such that &\ $\exists$ $C\geq 0$ such that \\
 $\log||\sigma\cdot w||^2\geq -C$   &\ $\log {||\sigma\cdot w||^2}-\log{||\sigma\cdot v||^2}  \geq -C $\\  
  all $\sigma\in G$ &\ all $\sigma\in G$ \\ 
  \hline \rule{0pt}{3ex}
  $G\cdot w$ \emph{closed} and $G_{w}$ finite & $\exists m\in\mathbb{N}$ such that $({I}^{q}\otimes v^m,w^{m+1})$ is semistable \\
 \hline
  \end{tabular} 
 \end{table}
  
   %%%%%%%%%%%%%%%%%%%%%%%%%%%%%%%%%%%%%%%%%%%%%%%%%%%%%%%%%%%%%
\section{Stability of Projective Varieties}
Fix $L\subset \mathbb{C}^{N+1}$ , $\dim(L)=n+1<N+1$.  Choose  $l\in \mathbb{N}$ satisfying $0\leq l\leq n$. Consider the Zariski open subset $\mathscr{U}_L$ of the Grassmannian defined by
\begin{align}
\begin{split}
 \mathscr{U}_L:=\{E\in G(N-l \ ,\ \mathbb{C}^{N+1})\ |\  H^{\bull}\left(0\ra E\cap L\ra E \overset{\pi_{L}}{\ra}\mathbb{C}^{N+1}/L\ra 0 \right)=0\}\ .
 \end{split}
 \end{align}
Observe that  $E\in\mathscr{U}_L$ if and only if
\begin{align}
\dim(\pi_{L}(E))=N-n \ .
\end{align}
Consider the subvariety $Z_L$  defined by
\begin{align}
Z_L:=G(N-l \ ,\ \mathbb{C}^{N+1})\setminus \mathscr{U}_L  \ .
 \end{align}
Then $E\in Z_L$ if and only if $-{\dim}(\pi_L(E))>n-N$ .

The rank plus nullity theorem implies that  
\begin{align}
\dim(E\cap L)+\dim(\pi_{L}(E))=N-l \ 
\end{align}
for {any} $E\in G(N-l \ ,\ \mathbb{C}^{N+1})$ .

 Therefore $E\in Z_L$ if and only if 
 \begin{align}
 \dim(E\cap L)>N-l+n-N=n-l \ .
 \end{align} 
 Therefore
 \begin{align}
 Z_L=\{ E\in G(N-l \ ,\ \mathbb{C}^{N+1})\ |\dim(E\cap L)\geq n-l+1\}\ .
 \end{align}
 
 Now we apply the previous linear algebra to a projective variety $X^n\subset\cpn$. Recall that for any $p\in X$ that the \emph{embedded tangent space} to $X$ at $p$ is the $n$ dimensional \emph{projective} linear subspace 
 \begin{align}
 \mathbb{T}_p(X)\in \mathbb{G}(n\ ,\ \cpn) \ 
 \end{align}
  obtained (for example) by projectivizing the tangent space the the cone over $X$ at any point $v\in\mathbb{C}^{N+1}\setminus \{0\} $ lying over $p$. 
  
  Given any $0\leq l \leq n$ we define the following subvariety $Z_{l}(X)$ of the Grassmannian by
 \begin{align}
 Z_{l}(X):=
 \{E\in \mathbb{G}(N-(l+1),\cpn)\ | \  \exists \ p\in X\cap E\ \mbox{\emph{and}}
 \ \dim(E\cap \mathbb{T}_p(X))\geq n-l\} \ .
 \end{align}
 Generally $Z_{l}(X)$ has \emph{codimension one}  in $\mathbb{G}(N-(l+1),\cpn)$ .

 To make the {defining polynomial} of $Z_{l}(X)$ concrete we view the Grassmannian in Stiefel coordinates \cite{sturmfels2017} by observing that there is a dominant map \footnote{The superscript $o$ denotes matrices of maximal rank.}
 \begin{align}
 M_{(l+1)\times (N+1)}^{o}\ni A\ra \pi({\tt{ker}}(A))\in \mathbb{G}(N-(l+1),\cpn)\ .
 \end{align}
 
 We may then consider the divisor (also denoted by $Z_{l}(X)$ )
 \begin{align}
 \overline{\pi^{-1}(Z_{l}(X))}\subset M_{(l+1)\times (N+1)} \ .
 \end{align}
 
Our ``new''  $Z_{l}(X)$ is now an irreducible algebraic hypersurface  in an affine space and hence is cut out by a single polynomial.
%$\Delta_{l+1}$
%\[ Z_{l}(X)=\{ (a_{ij}) \in  M_{(l+1)\times (N+1)}\ |\ \Delta_{l+1}(a_{ij})=0\}\ .\]

%\begin{remark} 
%\emph{We should point out that} \footnote{The height discrepancy arises naturally from this point of view.} \emph{$Z_{l+1}(X)$ is dominated by the variety of zeros of a larger system $I_\mathscr{X}$  in more variables $p\in \mathscr{X}$ where $\mathscr{X}$ is an auxiliary projective variety naturally associated to $X$.}
% \[ I_\mathscr{X}=\{(p,\ (a_{ij}))\in \mathscr{X}\times M_{(l+1)\times (N+1)}\ |\ s_k(p,(a_{ij}))=0 \ ; \ 1\leq k\leq m \} \]
%\emph{The situation can be visualized as follows}
%\begin{align} 
%\xymatrix{  
%	I_\mathscr{X} \ar@{^{(}->}[r]^-{\iota} \ar[d]^{{\pi_2}|_{I_\mathscr{X}}} &   {\mathscr{X}}\times M_{(l+1)\times (N+1)}\ar[d]^{\pi_2} \\  
%	Z_{l+1}(X)\ar@{^{(}->}[r] & M_{(l+1)\times (N+1)} \ .
%	 }
%\end{align}

%\emph{In geometric terms $(p,(a_{ij}))\in I_\mathscr{X}$ if and only if ${\tt{ker}}(a_{ij})$ fails to meet $\mathscr{X}$ generically at $p$ (and possibly at some other point $q$ ) .
% $Z_{l+1}(X)$ is therefore the \emph{resultant system} obtained by eliminating the variable $p$ from $I_{\mathscr{X}}$. Since $\mathscr{X}$ is projective, $Z_{l+1}(X)$ is a subvariety of $M_{(l+1)\times (N+1)}$ (see \cite{gkz} , \cite{weyman}, \cite{paulsergiou2019} ) .}
% \end{remark}

 %%%%%%%%%%%%%%%%%%%%%%%%%%%%%%%%%%%%%%%%%%%%%%%%%%%%%%%%%%%%%
 \subsection{Resultants}
Let $X^n\subset\cpn$ be an irreducible, $n$-dimensional, linearly normal, complex projective variety of degree $d$ .
  \begin{definition} (Cayley 1840's) \emph{The \textbf{\emph{associated hypersurface}} to $X^n\subset \cpn$ is given by}
\begin{align}
Z_{n}(X)=\{L\in \mathbb{G}(N-n-1,N)\ | L\cap X\neq \emptyset \} \ .
\end{align}
\end{definition}
As we have remarked, it is known that $Z_{n}(X)$ enjoys the following properties \\
\ \\
$i)$   $Z_{n}(X)$ is a {divisor} in $\mathbb{G}(N-n-1,N)$ ( and hence $M_{(n+1)\times (N+1)}$ ) . \\
 \ \\
$ii)$   $Z_{n}(X)$ is irreducible . \\
  \ \\
$iii)$  $\deg(Z_{n}(X))=d$ ( $=d(n+1)$ in Steifel coordinates ) . \\
  \ \\
  Therefore there exists $R_X\in H^0(\mathbb{G}(N-n-1,N), \mathcal{O}(d))$ such that
 \begin{align}
 \{R_X=0\}=Z_{n}(X)\ .
 \end{align}
 $R_X$ is the Cayley-{Chow} form of $X$. Modulo scaling,  $R_X$ is unique . Following the terminology of Gelfan'd \cite{gkz} we call $R_X$ the \textbf{$X$-\emph{resultant}} .  We will always view $R_X$ as a polynomial \footnote {It is necessarialy invariant under the natural action of $SL(n+1,\mathbb{C})$ .} in the matrix entries   
\begin{align}
R_X\in \mathbb{C}_{d(n+1)}[M_{(n+1)\times (N+1)}]^{SL(n+1,\mathbb{C})} \ . 
\end{align}
%%%%%%%%%%%%%%%%%%%%%%%%%%%%%%%%%%%%%%%%%%%%%%%%%%%%%%% 
\subsection{Hyperdiscriminants }  
Assume that $X\subset\cpn$ has degree $d\geq 2$. Let $X^{sm}$ denote the smooth points of $X$. For $p\in X^{sm}$ let 
$\mathbb{T}_p(X)$ be the {embedded tangent space} to $X$ at $p$ .
\begin{definition}
 \emph{ The \emph{\textbf{dual variety}} of $X$, denoted by $X^{\vee}$, is the Zariski closure of the set of \emph{tangent hyperplanes} to $X$ at its smooth points }
\begin{align}
X^{\vee}=\overline{ \{ f\in {\cpn}^{\vee} \ |  \  \mathbb{T}_p(X)\subset {\tt{ker}}(f) \ , \ p\in X^{sm}\} } \ .
\end{align}
\end{definition}
 Generally $X^{\vee}$ is codimension one in $ {\cpn}^{\vee}$. This holds, for example, whenever $X$ is a (nonlinear) projective curve or surface. Observe that we have the identity $$ X^{\vee}=Z_1(X)\ . $$ 
 
 For the purposes of understanding the Mabuchi energy, what is important is not the dual variety $X^{\vee}$ but the variety $Z_{n-1}(X)$.  This divisor also has a simple geometric description.
 $$Z_{n-1}(X)=\{ L\in \mathbb{G}(N-n \ ,\ \cpn)\ |\ \#(L\cap X)\neq \deg(X)\} $$
 
 It is known that $Z_{n}(X)$ enjoys the following properties \\
\ \\
$i)$   $Z_{n-1}(X)$ is a {divisor} in $\mathbb{G}(N-n,N)$ ( and hence $M_{n\times (N+1)}$ ) . \\
 \ \\
$ii)$   $Z_{n-1}(X)$ is irreducible . \\
  \ \\
$iii)$  $\deg(Z_{n-1}(X))=n(n+1)d-d\mu$ in Steifel coordinates  . \\
  \ \\
  Therefore there exists $\Delta_X\in H^0(\mathbb{G}(N-n,N), \mathcal{O}((n+1)d-d\frac{\mu}{n}))$ such that
 \begin{align}
 \{\Delta_X=0\}=Z_{n-1}(X)
 \end{align}
   Modulo scaling,  $\Delta_X$ is unique. Inspired by the terminology of Gelfan'd  we call $\Delta_X$ the \textbf{$X$-\emph{hyperdiscriminant}}.  We will always consider $\Delta_X$ as a polynomial \footnote {It is necessarily invariant under the natural action of $SL(n,\mathbb{C})$ .} in the appropriate matrix entries   
\begin{align}
\Delta_X\in \mathbb{C}_{n(n+1)d-d\mu}[M_{n\times (N+1)}]^{SL(n,\mathbb{C})} \ . 
\end{align}
 
 A word on notation is appropriate here. In \cite{gkz} the symbol $\Delta_X$ is used to denote \emph{the $X$-discriminant}. That is, the defining polynomial (when it exists) of the \emph{dual variety} $Z_0(X)$ of $X\subset \cpn$, whereas in this article  $\Delta_X$ is used to denote the defining polynomial of $Z_{n-1}(X)$. In \cite{sturmfels2017}  the defining polynomial of $Z_{n-1}(X)$ is denoted by $Hu_{X}$ and is called the \emph{Hurwitz form of $X\subset\cpn$}. The hyperdiscriminant and the Hurwitz form are the same polynomial.
 
 We summarize our constructions:\\
 Let $X^n\subset\cpn$ be a smooth, linearly normal complex projective variety.   We may associate two divisors  $Z_n(X)$ and $Z_{n-1}(X)$ cut out by irreducible polynomials $R_X$ and $\Delta_X$ respectively
\begin{align}
 \begin{split}
 &R_X\in  \mathbb{C}_{d(n+1)}[M_{(n+1)\times (N+1)}]^{SL(n+1,\mathbb{C})} \\
 &\Delta_X \cong \mathbb{C}_{n(n+1)d-d\mu}[M_{n\times (N+1)}]^{SL(n,\mathbb{C})} \ .
 \end{split}
 \end{align}
  For our purpose we must  {normalize the degrees} of these polynomials. From this point on we are interested in the pair 
\begin{align}
( {R_X}^{\ \deg(\Delta_X)}\ ,\  {\Delta_X}^{\ \deg(R_X)} )\ .
\end{align}
 
 Now we are prepared to make the following definitions. $X$ will always denote a smooth, linearly normal subvariety of $\cpn$.
 \begin{definition} \label{semistability} 
 {$X$ is   {\textbf {semistable}}  if and only if the pair $(R_X^{\ \deg(\Delta_X)},\Delta_X^{\ \deg(R_X)})$  is semistable for the action of $G$. Explicitly, the orbit closures are disjoint }
 \begin{align}\label{orbits}
\overline{\mathcal{O}}_{R\Delta}\cap \overline{\mathcal{O}}_{R}=\emptyset \ .
\end{align}
 \end{definition}
In (\ref{orbits}) we have defined
\begin{align}
&\mathcal{O}_{R\Delta}:=G\cdot[(R_X^{\ \deg(\Delta_X)}\ ,\ \Delta_X^{\ \deg(R_X)})]\quad, \quad \mathcal{O}_{R}:=G\cdot[(R_X^{\ \deg(\Delta_X)}\ ,\ 0)]\ .
\end{align}

 \begin{definition}\label{stability}
 {$X$ is  {\textbf {stable}}  if and only if the pair $(R_X^{\ \deg(\Delta_X)},\Delta_X^{\ \deg(R_X)})$ is stable for the action of $G$.  Explicitly, there is an integer $m\geq 2$ such that the pair $$
 ( {I}^q\otimes R_X^{(m-1)\deg(\Delta_X)}\ , \ \Delta_X^{m\deg(R_X)})$$
 is semistable for the action of $G$ and $q=\deg(R_X)\deg(\Delta_X)$.}
 \end{definition}

Proposition \ref{finiteautomorphism} immediately implies the following corollary
\begin{corollary} \label{stablepairisfinite}
The automorphism group of a stable variety is finite.
\end{corollary}

The reader should note that (semi)stability is independent of which lifts of  $R$ or $\Delta$ are chosen. Recall that we can only construct the \emph{divisors}, there will always be a scalar ambiguity in the choice of defining polynomial.

%%%%%%%%%%%%%%%%%%%%%%%%%%%%%%%%%%%%%%%%%%%%%%%%%%%%%%%%%%%%%%

\section{Asymptotic (semi)stability of Polarized Varieties }
We begin this section with a brief discussion of the well known equivalence among the various of $L^p$ norms on  spaces of homogeneous polynomials.

Given a homogeneous degree $d$ polynomial $P$ on $\mathbb{C}^{n+1}$ we identify it with a section of $\mathcal{O}(d)$ over $\cpn$. If we fix a Hermitian metric on $\mathbb{C}^{n+1}$ recall that the pointwise norm $|P|_{h_{FS}^d}([z])$ is given by
\begin{align}
|P|^2_{h_{FS}^d}([z]):=\frac{|P(z_0,\dots,z_n)|^2}{(|z_0|^2+\dots |z_n|^2)^d} \ .
\end{align}

For any $p\in [0,\infty]$ we define the $L^p$ norms by
\begin{align}
\begin{split}
&||P||_{0}:=\exp\left(\int_{\cpn}\log|P|_{h_{FS}^d}\om_{FS}^N\right)\ ,\ ||P||_{\infty}:=\sup_{[z]\in\cpn}|P|_{h_{FS}^d}([z])\\
&||P||_p:=\left(\int_{\cpn}|P|^p_{h_{FS}^d}\om_{FS}^N\right)^\frac{1}{p}\qquad p\in (0,\infty) \ .
\end{split}
\end{align}
{\remark{These satisfy the triangle inequality only for $p\in [1,\infty]$ .}}

\ \\
 Observe that $\log||P||_0$ is the \textbf{\emph{logarithmic Mahler measure}} of $P$ .   \ \\

 The following proposition is well known, we provide a simple proof below.
\begin{proposition}\label{arestov}(see \cite{faltings} ,\cite{lelong},\cite{bostgilletsoule} )
For any homogeneous polynomial $P$ of degree $d$ on $\cpn$ we have
\begin{align}\label{arestovinequality}
-\frac{d}{2}\left(\sum_{j=1}^{N}\frac{1}{j}\right)+\log||P||_{\infty}\leq \log||P||_{0} \leq \log||P||_{\infty} \ .
\end{align}
  \end{proposition}
 \begin{proof} 
  The content of the inequality (\ref{arestovinequality}) is the left hand side since the sup norm clearly dominates any $L^p$ norm. Recall that the mean zero Green's function for the scalar Fubini-Study Laplacian on $\cpn$ is given by (where $\rho$ denotes the geodesic distance between two points )
\begin{align}
G_{g_{FS}}(\rho)=\frac{1}{2N}\left(\sum_{j=1}^{N-1}\frac{1}{(N-j)\sin^{2N-2j}(\rho)}-2\log\sin(\rho)+\frac{1}{N}-2\sum_{j=1}^{N}\frac{1}{j}\right)\ .
\end{align}
In particular $$ G_{g_{FS}}(\rho) \geq -\frac{1}{N}\sum_{j=1}^{N}\frac{1}{j} \ . $$
For any (homogeneous) polynomial $P$ the Green's representation formula gives  
\begin{align}
\log|P|^2_{h_{FS}^d}([w])=\int_{\cpn}\log|P|^2_{h_{FS}^d}\om^N_{FS}-\int_{\cpn}G_{g_{FS}}\Delta \log|P|^2_{h_{FS}^d}\om^N_{FS}\qquad \mbox{for $[w]\notin Z(P)$}\ .
\end{align}
Since $G_{g_{FS}}$ has mean zero we have 
\begin{align}
\int_{\cpn}G_{g_{FS}}\Delta \log|P|^2_{h_{FS}^d}\om^N_{FS}=\int_{\cpn}G_{g_{FS}}(dN+\Delta \log|P|^2_{h_{FS}^d})\om^N_{FS}\ .
\end{align}
Plurisubharmonicity of $\log|P|^2_{h_{FS}^d}$ gives
\begin{align}
dN+\Delta \log|P|^2_{h_{FS}^d}\geq 0 \ .
\end{align}
 Therefore
 \begin{align}
 G_{g_{FS}}(dN+\Delta \log|P|^2_{h_{FS}^d})\geq -\frac{1}{N}\sum_{j=1}^{N}\frac{1}{j}(dN+\Delta \log|P|^2_{h_{FS}^d})\ .
 \end{align}
 Integrating this inequality gives the result.
 \end{proof}
 As a consequence of Proposition \ref{arestov} we deduce the following.
\begin{corollary}
For all $p\in (0 , \infty)$ we have
 \begin{align}\label{LpL0}
-\frac{d}{2}\left(\sum_{j=1}^{N}\frac{1}{j}\right)+\log||P||_{p}\leq \log||P||_{0} \leq \log||P||_{p} \ .
\end{align}
\end{corollary}
\noindent The right hand side of (\ref{LpL0}) comes from Jensen's inequality.
 In particular we get that the $L^2$ norm and the Mahler measure are equivalent, which is all that we need.
 Also observe that the inequality in Proposition (\ref{arestov}) becomes an equality for $P(z_0,\dots,z_N)=z_0^d$. \\

Now we return to the situation of $X\subset \cpn$. We assume that $X$ is smooth and linearly normal. Let $R_X$ and $\Delta_X$ denote the resultant and (hyper)discriminant. We remind the reader that these polynomials are only given \emph{up to scale}.\
   Propositions \ref{arestov} and \ref{logdistance}  justify the following definitions.
  \begin{definition}
 Let $X\subset \cpn$ .  Let $p\in [0,\infty]$. Choose any $L^p$ normalized $R$ and $\Delta$ . Then the $L^p$ \textbf{distance between the points} 
 \begin{align}
 \sigma\cdot[(R_X^{\ \deg(\Delta_X)}\ ,\ \Delta_X^{\ \deg(R_X)})] \ \mbox{and}\ \sigma\cdot[(R_X^{\ \deg(\Delta_X)}\ ,\ 0)]
  \end{align}
  is defined by
  \begin{align}
\log\tan {\tt{dist}}_p(\sigma):=    \log ||\sigma\cdot \Delta_X^{\ \deg(R_X)}||_p - \log||\sigma\cdot R_X^{\ \deg(\Delta_X)}||_p \ .
\end{align}
 \end{definition}
  
   \begin{definition}
 The  $L^p$ \textbf{distance between the orbit closures } is defined to be
 \begin{align}
 \log\tan{\tt{dist}}_p\left(\overline{\mathcal{O}}_{R\Delta},\overline{\mathcal{O}}_{R}\right):=
 \inf_{\sigma\in G} \log\tan {\tt{dist}}_p(\sigma)\ .
 \end{align}
 \end{definition}
 
 The point is that \emph{all} of the $L^p$ distances measure the same thing: any one of them detects the semistability of $X\subset \cpn$. What is extraordinary is that the infimum of the Mabuchi energy restricted to the Bergman metrics at level $k$ is \emph{exactly} the distance between the orbit closures in the $L^0$ distance .
  
Now we are prepared to the introduce \emph{asymptotic} (semi)stability of a polarized manifold $(X,L)$. We require an auxiliary Hermitian metric $h$ on $L$ with positive curvature $\om_h$.  The definition of asymptotic (semi)stability is independent of which $h$ is chosen. We must scale $R_X$ and $\Delta_X$ to have unit length in the norm $||\cdot||_0$.
%%%%%%%%%%%%%%%%%%%%%%%%%%%%%%%%%%%%%%%%%%%%%%%%%%%%%%%%%%%
\begin{definition} \label{asymptoticsemistability} 
 \emph{A polarized manifold $(X, {L})$ is \emph{\textbf { asymptotically semistable}} if and only if there is a uniform constant $C=C(h)\geq 0$ such that} 
 \begin{align}
 {\tt{dist}}_0(\overline{\mathcal{O}}_{R\Delta},\overline{\mathcal{O}}_{R})\succsim \exp(-Cd^2) \ 
 \end{align}
 \emph{for all sufficiently large $ L^k$-embeddings of degree} $d=k^n$ .
 \end{definition}
The author's previous work \cite{paul2012} shows that the orbit closures \emph{must} be disjoint for all powers of $L$, otherwise the Mabuchi energy is unbounded from below and no canonical metric exists.  Asymptotic semistability not only requires orbit closure separation for each embedding, but also that  {the orbit closures are not allowed to approach one another too quickly in the Mahler metric as the degree of the embedding increases}.

 \begin{definition} \label{asymptoticstability} 
 {A polarized manifold $(X,L)$ is {\textbf {asymptotically stable}} if and only if there are uniform constants $m\in\mathbb{Z}_{>0}$ and $C=C(h,m)$ such that} 
 \begin{align}
 {\tt{dist}}_0(\overline{\mathcal{O}}_{(v,w)},\overline{\mathcal{O}}_ {v})\succsim \exp(-Ck^{2n+1}) \ 
 \end{align}
 {for all sufficiently large $k$ (the power of the embedding) .}
 $$(v,w):= ({I}^q\otimes R_X^{(km-1)\deg(\Delta_X)}\ , \ \Delta_X^{km\deg(R_X)}) \ .$$
 \end{definition}
 {As in the definition of asymptotic semistability, both $R_X$ and $\Delta_X$ have been scaled to have length one in the norm $||\cdot||_0$. The reader should observe that the speeds of approach of the orbit closures in the definitions and asymptotic stability and semistability differ by a single factor of $k$ . }
 
%%%%%%%%%%%%%%%%%%%%%%%%%%%%%%%%%%%%%%%%%%%%%%%%%%%%%%%%%%%%
 
%%%%%%%%%%%%%%%%%%%%%%%%%%%%%%%%%%%%%%%%%%%%%%%%%% 

\section{ Asymptotic Stability and Properness of the Mabuchi Energy}
 We recall some definitions surrounding Mabuchi's K-energy map.
Let 
\begin{align}  
(X^n,\om)  \ , \ n={\tt{dim}}_{\mathbb{C}}(X) \ 
\end{align}
be a compact K\"ahler manifold. Recall that the K\"ahler form $\om$ is given locally by a Hermitian positive definite matrix of functions
  \begin{align}
  \om=\frac{\sqrt{-1}}{2\pi}\sum_{ i,j }g_{i\overline{j}}dz_{i}\wedge d\overline{z}_j  \ .
  \end{align}
The Ricci form of $\om$ is the smooth $(1,1)$ form on $X$ given by
\begin{align}
\begin{split}
 {\tt{Ric}}(\om):= \frac{-\sqrt{-1}}{2\pi }\dl\dlb\log\det(g_{i\overline{j}}) =\sum_{i,j}\frac{-\sqrt{-1}}{2\pi }R_{i\overline{j}}dz_{i}\wedge d\overline{z}_j  \ .
\end{split}
\end{align}
The scalar curvature is by definition the contraction of the Ricci curvature 
 \begin{align}
{\tt{Scal}}(\om): =\sum_{i,j}g^{i\overline{j}}R_{i\overline{j}}\in C^{\infty}(X)\ .
 \end{align}  
The volume $V$ and the average of the scalar curvature $\mu$ depend only on $[\om]$ and are given by
\begin{align}
 V=\int_X\om^n \ , \ \mu=\frac{1}{V}\int_X{\tt{Scal}}(\om)\om^n   \ .
 \end{align}
The space of K\"ahler metrics in the class $[\om]$ is defined by
\begin{align}
\begin{split}
& \mathcal{H}_{\om}:=\{\vp\in C^{\infty}(X)\ |\ \om_{\vp}:=\om+\frac{\sqrt{-1}}{2\pi}\dl\dlb\vp  >0 \}  \ .
\end{split}
\end{align} 
  \begin{definition} (Mabuchi \cite{mabuchi}) 
 \emph{The \textbf{\emph{K-energy map}} $\nu_{\om}:\mathcal{H}_{\om}\ra \mathbb{R}$ is given by
\begin{align}
\nu_{\om}(\varphi):=-\frac{1}{V}\int_0^1\int_{X}\dot{\varphi_{t}}({\tt{Scal}}(\om_{\varphi_{t}})-\mu)\omega_{t}^{n}dt 
\end{align}
 $\vp_{t}$ is a $C^1$ path in $\mathcal{H}_{\om}$ satisfying $\vp_0=0$ , $\vp_1=\vp$ .}
 \end{definition}
Mabuchi shows that  $\nu_{\om}$ is independent of the path chosen. It is clear that 
 {$\vp$ is a critical point for $\nu_{\om}$ if and only if $$ {\tt{Scal}}(\om_{\varphi})\equiv \mu\ . $$ }
 What is relevant for the present article is the following theorem, first established by Bando and Mabuchi in the case $L=-K_{X}$ , and then generalized some years later by Donaldson and Li .   
 
 \begin{theorem} (see \cite{bando-mabuchi87}, \cite{donaldson2001}, \cite{donaldson2005}, \cite{chili2011} ) Let $(X,L)$ be a polarized manifold, and assume that there is a constant scalar curvature metric in the class ${c_1}(L )$. Then the Mabuchi energy is bounded below on $\mathcal{H}_{\om}$ where $h$ is any Hermitian metric on $L$ with positive curvature $\om$.
\end{theorem}
 
We recall the Aubin $J_{\om}$ functional (see \cite{aubin84}) and the associated energy $F_{\om}^o$  
\begin{align}
\begin{split}
& J_{\omega}(\varphi):= \frac{1}{V}\int_{X}\sum_{i=0}^{n-1}\frac{\sqrt{-1}}{2\pi}\frac{i+1}{n+1}\dl\varphi \wedge \dlb
\varphi\wedge \omega^{i}\wedge {\omega_{\varphi} }^{n-i-1}\\
& F_{\om}^o(\varphi):=J_{\omega}(\varphi)-\frac{1}{V}\int_{X}\varphi\ \om^n \ .
\end{split}
\end{align}

\begin{definition}\label{proper}(Tian \cite{tian97})
{Let $(X,\om) $ be a K\"ahler manifold . The Mabuchi energy is {\textbf{proper}} provided there exists constants $a >0$ and $b$ such that for all $\varphi\in \mathcal{H}_\om$ we have }
\begin{align}
\begin{split}
&\nu_{\om}(\vp)\geq a J_{\omega}(\varphi)+b  \ .
 \end{split}
\end{align}
\end{definition}

Let $(X ,L)$ be a polarized manifold. Let $h$ be a smooth Hermitian metric on $L$ with positive curvature $\om$. Choose $k$ large enough so that there is an embedding
\begin{align}
\iota_k:X\ra \mathbb{P}(H^0(X , L^k)^*) \ .
\end{align}
We will always assume that the embedding is given by a \emph{unitary} basis of sections $\{S_i \}$. Similarly we outfit $H^0(X , L^k)$ with the Hodge $L^2$ inner product. We let $\om_{FS}$ denote the corresponding Fubini-Study K\"ahler metric on the (dual) projective space of sections. Then
\begin{align}
\iota^*_{k}\om_{FS}|_{\iota_k(X)}=k\om_h+\frac{\sqrt{-1}}{2\pi}\dl\dlb\log\left(\sum^{N_k}_{i=0}|S_i|^2 \right)\ .
\end{align}
Let $G=SL(H^0(X , L^k))$ , then $\sigma\in G$ acts on the sections by
\begin{align}
\sigma\cdot S_{i}=\sum_{0\leq j\leq N_k}\sigma_{ij}S_j \ .
\end{align}
Define 
\begin{align}
\Psi_{\sigma}:=\log \sum_{0\leq i\leq N_k+1}|\sigma\cdot S_i|^2\ .
\end{align}
The \emph{Bergman metrics} of level $k$ are given by
\begin{align}
\mathscr{B}_{N_k}:=\{ \frac{1}{k}\Psi_{\sigma} \ |\ \sigma\in SL(N_k+1,\mathbb{C}) \} \subset \mathcal{H}_{\om}\ .
\end{align}
In the discussion below we need to distinguish between the potential $\Psi_{\sigma}$ which is a function on $X$, and the closely related potential $\vps$, which is a function on $\mathbb{P}^{N_k}$.  Note that $\vps$ is  a K\"ahler potential on $\iota_k(X)$ relative to the restriction $\om_{FS}|_{\iota_k(X)}$. These two potentials are related as follows. 
\begin{align}
\begin{split}
&\sigma^{*}\om_{FS}=\om_{FS}+\frac{\sqrt{-1}}{2\pi}\dl\dlb\vps \ , \ \vps([z])=\log\frac{|\sigma\cdot z|^2}{|z|^2} \\
&\iota^*_k(\sigma^{*}\om_{FS}|_{\iota_k(X)})=k\om_h+\frac{\sqrt{-1}}{2\pi}\dl\dlb \Psi_{\sigma} \\
&\Psi_{\sigma}=\vps\circ \iota_k+\log \sum_{0\leq i\leq N_k+1}|S_i|^2 \ .
\end{split}
\end{align}
 
 A key ingredient in this paper is the following result of Tian \cite{tianberg} .
 \begin{theorem} (Tian's Thesis)
 The spaces $\mathscr{B}_{N_k}$ are dense in the $\tt{C}^2$ topology
 \begin{align}
 \overline{\bigcup _{k}\mathscr{B}_{N_k}}=\mathcal{H}_{\om} \ .
 \end{align}
\end{theorem}

Now we are prepared to establish that the asymptotic (semi)stability of $(X,L)$ is equivalent to the (lower bound) global coercive estimate for the Mabuchi energy $\nu_{{\om}_{h}}$  for any ${\om}_{h}\in{c}_1(L)$ .
  We need to compare the Mabuchi and Aubin energies of the reference metric $\om:=\om_{h}$ with the restrictions of the Fubini-Study metrics coming from the large projective embeddings. It is easy to see that $\nu_{\om}$ does not scale but $F_{\om}^o$ does scale as we pass between $\om$ and ${\om_{FS}|_{\iota_k(X)}}$. We collect the precise comparisons below, where $o(1)$ denotes any quantity that converges to 0 as $k\ra \infty$. The $o(1)$'s below have the form $O(\frac{\log(k)}{k})$. 
  
 \begin{align}\label{comparisons}
 \begin{split}
 & \nu_{\om}\left(\frac{\Psi_{\sigma}}{k}\right)=\nu_{\om_{FS}|_{\iota_k(X)}}\left(\varphi_{\sigma}\right)+o(1) \ , \  J_{\om}\left(\frac{\Psi_{\sigma}}{k}\right)=\frac{1}{k}J_{\om_{FS}|_{\iota_k(X)}}\left(\varphi_{\sigma}\right)+o(1) \\
 & \int_{X}\frac{\Psi_{\sigma}}{k}\frac{\om^n}{V_o}=\frac{1}{V}\int_{\iota_k(X)}\frac{\varphi_{\sigma}}{k}\om^n_{FS}+o(1) \ .
 \end{split}
 \end{align}

All of the results in this article depend on the following theorem which completely describes the Mabuchi energy restricted to the space of Bergman metrics associated to the embedding $X^n\subset\cpn$.   

 \ \\
 \textbf{Theorem A .} (\cite{paul2012})
\emph{ 
 There is a norm  $||\cdot ||$ on the space of polynomials such that }
\begin{align} \label{norm}
\begin{split}
& \ d^2(n+1)\nu_{\om_{FS}|_X}(\varphi_{\sigma})= {\deg(R_X)}\log  \frac{{||\sigma\cdot\Delta_{X}||}^{2}}{{||\Delta_{X}||}^{2}} - {\deg(\Delta_{X})} \log  \frac{{||\sigma\cdot R_{X}||}^{2}}{||R_{X}||^2} \ .
\ \\
\end{split}
\end{align}
 The norm appearing in (\ref{norm}) was first considered by Gang Tian in his early works on CM stability \cite{bottchrnfrms}, \cite{tian94}, \cite{tian97}, \cite{tiansurv}, \cite{tianbook} . This norm is conformally equivalent to the $L^2$ norm with a continuous potential $\theta$
 \begin{align}
 || \cdot ||:= e^{\theta}||\cdot ||_{L^2}\ .
 \end{align}
 In the situation considered by Tian in \cite{tian97} it seems there is little one could say about $\theta$  beyond it's (H\"older) continuity. However, for families of \emph{divisors}, the situation considered here,  $\theta$ can be described explicitly which allows us to significantly improve Theorem A.  
 
\begin{theorem}\label{precisetheorem}
 {Let $X\subset \cpn$ be a smooth, linearly normal complex projective variety then the following holds for all $\sigma\in G$}
 \begin{align}  \label{precisetheoremA}
\begin{split}
& \ d^2(n+1)\nu_{\om_{FS}|_X}(\varphi_{\sigma})= {\deg(R_X )}\log   \frac{{||\sigma\cdot\Delta_X||_0}^{2}}{||\Delta_X||_0^2}  - {\deg(\Delta_X)} \log   \frac{{||\sigma\cdot R_X||_0}^{2}}{||R_X||_0^2}  \ .
\ \\
\end{split}
\end{align}
 \end{theorem}
 
\begin{proof}
  We identify the conformal factor $\theta$. Consider $\mathbb{C}^{n+1}$ equipped with it's standard metric. Let $d$ be any positive integer.  We identify the space of homogeneous polynomials of degree $d$ on $\mathbb{C}^{n+2}$ with
 $H^{o}(\mathbb{P}^{n},\mathscr{O}(d))$.   We let $B$ denote the corresponding complete linear system and $\mathscr{X}_d$ the universal family of hypersurfaces over $B$
 \begin{align}
 B:= \mathbb{P}( H^{o}(\mathbb{P}^{n},\mathscr{O}(d)))\quad \mathscr{X}_d:=\{([S] ,[w])\in B\times \mathbb{P}(W) \ | S(w)=0\ \} \ .
 \end{align}
  Using the projections
 \begin{align}   
    \xymatrix{
\mathscr{X}_d\ar[d]^{p_1} \ar[r]^{p_2}&\mathbb{P}^{n} \\
B          &}
 \end{align} 
 we define a closed positive $(1,1)$ current $u$ on $B$  by 
 $$u:={p_{1}}_{*}p_{2}^{*}(\om_{FS}^{n}) \ . $$
 Explicitly, for any smooth form $\alpha$ on $B$ of correct type we define
\begin{align}\label{current}
\int_{B}u\wedge \alpha:=\int_{ \mathscr{X}_d}p_2^*(\om^{n})\wedge p_1^*(\alpha) \ .
\end{align}
  Since DeRham and current cohomology on $B$ coincide we see at once that $[u]=[\om_B]$ where $\om_{B}$ is a smooth $(1,1)$ form on $B$.  

\begin{proposition}  (\cite{tian97} Lemma 8.7 pg. 32) 
There is a continuous function $\theta$ on $B$ such that, in the sense of currents we have
\begin{align}
u=\om_{B}+\frac{\sqrt{-1}}{2\pi}\dl\dlb \theta \quad \quad [\om_{B}]=[{c}_1(\mathscr{O}_B(1))] \ .
\end{align}
\end{proposition}
 With this said the conformal factor that appears on the right hand side of (\ref{norm}) is $e^{\theta}$ and for any section $S$ of $\mathscr{O}_{\mathbb{P}^{n+1}}(d)$ it's norm is defined to be
  \begin{align}\label{tiannorm}
 ||S||:= e^{\theta([S])}||S||_{L^2} \ .
 \end{align}
  This introduces a bounded ``error''  on the right hand side of (\ref{norm}) when we relate the Mabuchi energy to the $L^2$ norm. An explicit description of $\theta$ is obtained by noting that $\mathscr{X}_d$ is a divisor in $B\times\mathbb{P}^{n+1}$ cut out by a section $\Psi$ of $p_1^{*}\mathscr{O}_{B}(1)\otimes p_2^{*}\mathscr{O}_{\mathbb{P}^{n}}(d)$
   
\begin{align} 
\begin{split}
&\xymatrix{&p_1^{*}\mathscr{O}_{B}(1)\otimes p_2^{*}\mathscr{O}_{\mathbb{P}^{n+1}}(d)\ar[d]\\  
	\mathscr{X}_d \ar@{^{(}->}[r]^-{\iota}  &  B\times\mathbb{P}^{n} \ar@/^1pc/[u]^{ {\Psi}} } \\
\ \\
&	  \Psi ([S], [w]):=1_{S}\otimes S([w]) 
\end{split}	
	\end{align}
 Observe that, in the natural Hermitian metric on $p_1^{*}\mathscr{O}_{B}(1)\otimes p_2^{*}\mathscr{O}_{\mathbb{P}^{n}}(d)$ the $\log$ of the length of $\Psi$ is  
  \begin{align}
   \log|\Psi([S], [w])|_h^2=\log\frac{|S([w])|^2_{h^d_{FS}}}{||S||^2_{L^2}}\ .
  \end{align}
  Next observe that the Poincar\'e-Lelong formula gives
\begin{align}
\begin{split}
&\int_{ \mathscr{X}_d}p_2^*(\om^{n+1})\wedge p_1^*(\alpha)\\
&=\int_{B\times\mathbb{P}^{n}}\left(dp_1^*(\om_{FS})+p_2^*(\om_B)+\frac{\sqrt{-1}}{2\pi}\partial\overline{\partial}\log|\Psi|^2_{h}\right)\wedge p_2^*(\om_{FS}^{n})\wedge p_1^*(\alpha) \\
&=\int_B\left(\om_B+\frac{\sqrt{-1}}{2\pi}\partial_B\overline{\partial}_B\int_{\mathbb{P}^{n}}\log|\Psi|_h^2\om _{FS}^{n}\right)\wedge \alpha \ .
\end{split}
\end{align}
Therefore we get
\begin{align}\label{conformalfactor}
 \theta([S])=\int_{\mathbb{P}^{n}}\log\frac{|S([w])|^2_{h^d_{FS}}}{||S||^2_{L^2}}\om_{FS}^{n}\ .
 \end{align}
     Inserting (\ref{conformalfactor})  into the right hand side of (\ref{tiannorm}) gives that the norm of the discriminant  and resultant are given explicitly by
 \begin{align}
\begin{split}
&\log  \frac{{||\sigma\cdot\Delta||}^{2}}{{||\Delta||}^{2}}=\int_{\mathbb{P}(M_{n\times (N+1)})}\log|\sigma\cdot\Delta|^2_{h_{FS}} - \int_{\mathbb{P}(M_{n\times (N+1)})}\log|\Delta|^2_{h_{FS}}\\
&\log  \frac{{||\sigma\cdot R||}^{2}}{{||R||}^{2}}=\int_{\mathbb{P}(M_{n+1\times (N+1)})}\log|\sigma\cdot R|^2_{h_{FS}}  - \int_{\mathbb{P}(M_{n+1\times (N+1)})}\log|R|^2_{h_{FS}} \ .
 \end{split}
 \end{align}
In other words the norm $|| \cdot ||$ in (\ref{tiannorm}) is given by
\begin{align}
||S||=||S||_0 \quad \mbox{for any $S\in H^0(\mathbb{P}^{n},\mathscr{O}(d))$} \ .
\end{align}
\end{proof}
\noindent The second part of Theorem \ref{maintheorem}, namely equivalence between asymptotic semistability and a global lower bound for the Mabuchi energy, follows from (\ref{precisetheoremA}) , \cite{tianberg} and the corollary below.
  \begin{corollary}
 For any polarized manifold $(X,L)$ and any large $k$ embedding $ X\subset \cpn$ the infimum of the Mabuchi energy restricted to 
 $G=SL(N_k+1,\mathbb{C}) $ is given by
 \begin{align}\label{distanceformula}
 \inf_{{\sigma}\in G} d^2(n+1)\nu_{\om_{FS}|_{\iota_k(X)}}\left(\varphi_{\sigma}\right)= \log\tan{\tt{dist}}_{0}(\overline{\mathcal{O}}_{R\Delta} , \overline{\mathcal{O}}_{R} )\ .
 \end{align}
 \end{corollary}
\begin{proof}
This follows at once from the definition of the distance in the $L^p$ metrics.
\end{proof}
 \begin{remark}
 The reader should compare (\ref{distanceformula}) with the corollary on pg. 257 of \cite{paul2012} .
 \end{remark}
 Now we can show the first part of Theorem \ref{maintheorem}, namely equivalence between asymptotic stability and the global coercive estimate for the Mabuchi energy. 
\begin{proposition}
Let $m$ be a positive integer. For any polarized manifold $(X,L)$ and any large $k$ embedding we have
 \begin{align}\label{coercive}
\begin{split}
 \inf_{\frac{\Psi_{\sigma}}{k}\in\mathscr{B}_{N_k}} &\left( m\nu_{\om_{h}}\left(\frac{\Psi_{\sigma}}{k}\right)- 
\frac{\deg(\Delta_X)}{d} J_{\om_{h}}\left(\frac{\Psi_{\sigma}}{k}\right)\right) \\
&=\frac{ k^{-(2n+1)}}{(n+1)}\log\tan{\tt{dist}}_{0}(\overline{\mathcal{O}}_{(v,w)} , \overline{\mathcal{O}}_{v} )+O(1)\ ,
\end{split}
 \end{align}
 { where we have defined the pair $(v,w)$ and $q$ as follows} 
 $$ (v,w):= ( {I}^q\otimes R_X^{(km-1)\deg(\Delta_X)}\ , \ \Delta_X^{km\deg(R_X)})\ , \ q:=\deg(R_X)\deg(\Delta_X) $$
 and the distance ${\tt{dist}}_{0}$ in (\ref{coercive}) has been extended to $I$ by simply using the Hilbert-Schmidt norm on matrices.
 \end{proposition}

\begin{proof}
We begin with the following crucial observation, which was shown to the author by Gang Tian.
\begin{lemma}\label{tian'slemma}
There is a uniform constant $C$ such that for all sufficiently large $k\in \mathbb{N}$ we have
\begin{align}
C+\frac{1}{k}\log\left(\frac{||\sigma||^2}{N_k+1} \right)\leq \int_X\frac{\Psi_{\sigma}}{k}\frac{\om^n}{V_o} \quad  .
\end{align}
\end{lemma}
\begin{proof} 
If $||\sigma||^2:={\tt{Trace}}(\sigma\sigma^*)$  then we observe that the unitarity of the basis gives
\begin{align}
 \sum_{0\leq i\leq N_k}\frac{||\sigma\cdot S_{i}||^2}{||\sigma||^2}=1 \ .
\end{align}
Therefore there is an index $j$ such that
\begin{align}\label{indexinequality}
 \log||\sigma\cdot S_{j}||^2\geq \log \frac{||\sigma||^2}{N_k+1}\ .
\end{align}
Define $$T_j^{\sigma}:=\frac{\sigma\cdot S_{j}}{||\sigma\cdot S_{j}||} \ .$$ Let $\alpha(L)$ be Tian's alpha invariant \cite{tianferm} , and choose any $0<\beta < \alpha(L)$ then there exists a uniform constant $C(\beta)>0$ such that
\begin{align}
\int_X\left(\frac{1}{|T_j^{\sigma}|^2}\right)^{\frac{\beta}{k}}\frac{\om^n}{V} \leq C(\beta) \ .
\end{align}
Jensen's inequality gives 
\begin{align}
-\frac{\beta}{k}\int_{X}\left(\log|\sigma\cdot S_{j}|^2-\log||\sigma\cdot S_{j}||^2\right)\frac{\om^n}{V} \leq \log C(\beta) \ .
\end{align}
Equivalently
\begin{align}
\frac{\beta}{k} \log||\sigma\cdot S_{j}||^2 \leq \frac{\beta}{k}\int_{X}\log|\sigma\cdot S_{j}|^2+\log C(\beta) \ .
\end{align}
Applying inequality (\ref{indexinequality}) we see that
\begin{align}
\begin{split}\label{beta}
-\frac{1}{\beta}\log  C(\beta)+ \frac{1}{k}\log\left(\frac{||\sigma||^2}{N_k+1} \right)  &\leq \frac{1}{k}\int_X\log|\sigma\cdot S_j|^2\frac{\om^n}{V_o}\leq \int_{X}\frac{\Psi_{\sigma}}{k}\frac{\om^n}{V_o} \ .
 \end{split}
\end{align}
 \end{proof}
\noindent The comparison formulas (\ref{comparisons}) and the preceding lemma imply that
 \begin{align}\label{Jlognorm}
  J_{\om_{h}}\left(\frac{\Psi_{\sigma}}{k}\right)=\frac{1}{k}F^o_{\om_{FS}|_{\iota_k(X)}}(\vps)+\frac{1}{k}\log||\sigma||^2 +O(1) \ .
 \end{align}
Recall the well known proposition .\footnote{This amounts to expressing the {Faltings height} of $X$ in terms of the Mahler measure of the Cayley (Chow) form of $X$. This was first shown by P. Phillipon and independently by C. Soul\'e around 1991. See also \cite{zhang}, \cite{luo}, \cite{gacms} .}
\begin{proposition}\label{phillipon} (\cite{phillipon1991}, \cite{soule1991})
For any linearly normal projective variety $X\subset \cpn$ we have 
\begin{align}
-\deg(R_X)F^{o}_{\om_{FS}|_{\iota_k(X)}}(\varphi_{\sigma})=\log||\sigma\cdot R_X||_0 \ .
\end{align}
\end{proposition}
\noindent In the above Proposition we have chosen $R_X$ to have length one in the Mahler norm. Inserting (\ref{Jlognorm}) into Proposition \ref{phillipon} allows us to express $J_{\om_h}|_{\mathscr{B}_{N_k}}$ as a distance function
\begin{align}\label{preciseJ}
\begin{split}
&\frac{\deg(\Delta_X)}{d}J_{\om_{h}}\left(\frac{\Psi_{\sigma}}{k}\right)=\\
&\frac{1}{k^{2n+1}(n+1)}\left(-\deg(\Delta_X)\log||\sigma\cdot R_X||^2_0+q\log||\sigma||^2\right) +O(1) \ .
\end{split}
\end{align}
Theorem \ref{precisetheorem} and the comparison formulas (\ref{comparisons}) give
\begin{align}\label{precisemabuchi}
\begin{split}
&m \nu_{\om_h}\left(\frac{\Psi_{\sigma}}{k}\right)=\\
&\frac{1}{k^{2n+1}(n+1)}\left(km \deg(R_X)\log||\sigma\cdot\Delta_X||_0^2-km\deg(\Delta_X)\log||\sigma\cdot R_X||_0^2\right)+o(1) \ .
\end{split}
\end{align}
As usual, we have chosen representatives satisfying $||R_X||_0=||\Delta_X||_0=1$. 
Now subtract (\ref{preciseJ}) from (\ref{precisemabuchi}) and use the definition of the $L^0$ distance to get
\begin{align}\label{almostthere}
\begin{split}
&  m\nu_{\om_{h}}\left(\frac{\Psi_{\sigma}}{k}\right)-\frac{\deg(\Delta_X)}{d}J_{\om_{h}}\left(\frac{\Psi_{\sigma}}{k}\right) = \\
&\frac{k^{-(2n+1)}}{(n+1)}\log\tan{\tt{dist}}_0(\sigma\cdot[(v,w)],\sigma\cdot [(v,0)])+O(1) \ . \\
\end{split}
\end{align}
Recall that the pair $(v,w)$ is given by
\begin{align}
(v,w):= (I^q\otimes R_X^{(km-1)\deg(\Delta_X)},\Delta_X^{km\deg(R_X)}) \ .
\end{align}
 Taking the $\inf$ over $G$ on both sides of (\ref{almostthere}) completes the proof of Theorem \ref{maintheorem}.
 \end{proof}

\begin{center}\textbf{Acknowledgements}\end{center}
  The ideas in this paper are a direct outgrowth of Gang Tian's many important contributions to K\"ahler Geometry, and in particular his seminal, but unfortunately overlooked, concept of \emph{CM Stability}.  The author appreciates Tian's many patient explanations of his idea that the connection between algebraic geometry and canonical metrics should come from the Mabuchi energy restricted to the Bergman metrics. The author also thanks Professor Jeff Viaclovsky for a very careful reading of the article and for many suggestions which improved the exposition. The author also benefited from several conversations with Dr. Chi Li at Purdue in the Fall of 2019.
  
  This work was carried out during two visits to the Mathematics Department of the Massachusetts Institute of Technology. The first visit was in the spring of 2013 and the second in the spring of 2021, where this work was finally completed. The author thanks his host, Professor Tomasz Mrowka, for the kind invitations and for many useful discussions on the topic of the paper and his encouragement.     
\bibliographystyle{plain} 
\bibliography{ref.bib}

\begin{thebibliography}{10}

\bibitem{aubin84}
Thierry Aubin.
\newblock R\'{e}duction du cas positif de l'\'{e}quation de {M}onge-{A}mp\`ere
  sur les vari\'{e}t\'{e}s k\"{a}hl\'{e}riennes compactes \`a la
  d\'{e}monstration d'une in\'{e}galit\'{e}.
\newblock {\em J. Funct. Anal.}, 57(2):143--153, 1984.

\bibitem{bando-mabuchi87}
Shigetoshi Bando and Toshiki Mabuchi.
\newblock Uniqueness of {E}instein {K}\"ahler metrics modulo connected group
  actions.
\newblock In {\em Algebraic geometry, Sendai, 1985}, volume~10 of {\em Adv.
  Stud. Pure Math.}, pages 11--40. North-Holland, Amsterdam, 1987.

\bibitem{birkes71}
David Birkes.
\newblock Orbits of linear algebraic groups.
\newblock {\em Ann. of Math. (2)}, 93:459--475, 1971.

\bibitem{bgs1}
J.-M. Bismut, H.~Gillet, and C.~Soul{\'e}.
\newblock Analytic torsion and holomorphic determinant bundles. {I}.
  {B}ott-{C}hern forms and analytic torsion.
\newblock {\em Comm. Math. Phys.}, 115(1):49--78, 1988.

\bibitem{bgs2}
Jean-Michel Bismut, Henri Gillet, and Christophe Soul{\'e}.
\newblock Analytic torsion and holomorphic determinant bundles. {II}. {D}irect
  images and {B}ott-{C}hern forms.
\newblock {\em Comm. Math. Phys.}, 115(1):79--126, 1988.

\bibitem{bgs3}
Jean-Michel Bismut, Henri Gillet, and Christophe Soul{\'e}.
\newblock Analytic torsion and holomorphic determinant bundles. {III}.
  {Q}uillen metrics on holomorphic determinants.
\newblock {\em Comm. Math. Phys.}, 115(2):301--351, 1988.

\bibitem{bostgilletsoule}
J.-B. Bost, H.~Gillet, and C.~Soul\'{e}.
\newblock Heights of projective varieties and positive {G}reen forms.
\newblock {\em J. Amer. Math. Soc.}, 7(4):903--1027, 1994.

\bibitem{cds3}
Xiuxiong Chen, Simon Donaldson, and Song Sun.
\newblock K\"{a}hler-{E}instein metrics on {F}ano manifolds. {I}:
  {A}pproximation of metrics with cone singularities.
\newblock {\em J. Amer. Math. Soc.}, 28(1):183--197, 2015.

\bibitem{cds2}
Xiuxiong Chen, Simon Donaldson, and Song Sun.
\newblock K\"{a}hler-{E}instein metrics on {F}ano manifolds. {II}: {L}imits
  with cone angle less than {$2\pi$}.
\newblock {\em J. Amer. Math. Soc.}, 28(1):199--234, 2015.

\bibitem{cds1}
Xiuxiong Chen, Simon Donaldson, and Song Sun.
\newblock K\"{a}hler-{E}instein metrics on {F}ano manifolds. {III}: {L}imits as
  cone angle approaches {$2\pi$} and completion of the main proof.
\newblock {\em J. Amer. Math. Soc.}, 28(1):235--278, 2015.

\bibitem{chencheng1}
X.X. Chen and J.~Cheng.
\newblock On the constant scalar curvature {K}\"ahler metrics, apriori
  estimates.
\newblock {\em arXiv:1712.06697}, 2017.

\bibitem{chencheng2}
X.X. Chen and J.~Cheng.
\newblock On the constant scalar curvature {K}\"ahler metrics {II} existence
  results.
\newblock {\em arXiv:1801.00656}, 2018.

\bibitem{chencheng3}
X.X. Chen and J.~Cheng.
\newblock On the constant scalar curvature {K}\"ahler metrics {III}-general
  automorphism group.
\newblock {\em arXiv:1801.05907}, 2018.

\bibitem{cupit2003}
St{\'e}phanie Cupit-Foutou.
\newblock Classification of two-orbit varieties.
\newblock {\em Comment. Math. Helv.}, 78(2):245--265, 2003.

\bibitem{donaldson2005}
S.~K. Donaldson.
\newblock Scalar curvature and projective embeddings {II}.
\newblock {\em Q. J. Math.}, 56(3):345--356, 2005.

\bibitem{donaldson2001}
S.K. Donaldson.
\newblock Scalar curvature and projective embeddings {I}.
\newblock {\em JDG}, 59:479--522, 2001.

\bibitem{faltings}
G.~Faltings.
\newblock Diophantine approximation on abelian varieties.
\newblock {\em Annals of Mathematics}, 133:549--576, 1991.

\bibitem{gkz}
I.~M. Gelfand, M.~M. Kapranov, and A.~V. Zelevinsky.
\newblock {\em Discriminants, resultants, and multidimensional determinants}.
\newblock Mathematics: Theory \& Applications. Birkh\"auser Boston Inc.,
  Boston, MA, 1994.

\bibitem{lelong}
P.~Lelong.
\newblock Mesure de {M}ahler et calcul de constantes universelles pour les
  polynomes de n variables.
\newblock {\em Mathematische Annalen}, 299:673 --695, 1994.

\bibitem{chili2011}
Chi Li.
\newblock Constant scalar curvature {K}\"ahler metric obtains the minimum of
  {K}-energy.
\newblock {\em IMRN}, 9:2161 -- 2175, 2011.

\bibitem{luo}
Huazhang Luo.
\newblock Geometric criterion for {G}ieseker-{M}umford stability of polarized
  manifolds.
\newblock {\em J. Differential Geom.}, 49(3):577--599, 1998.

\bibitem{mabuchi}
Toshiki Mabuchi.
\newblock {$K$}-energy maps integrating {F}utaki invariants.
\newblock {\em Tohoku Math. J. (2)}, 38(4):575--593, 1986.

\bibitem{gacms}
S.~T. Paul.
\newblock Geometric analysis of {C}how {M}umford stability.
\newblock {\em Adv. Math.}, 182(2):333--356, 2004.

\bibitem{paul2012}
Sean~Timothy Paul.
\newblock Hyperdiscriminant polytopes, {C}how polytopes, and {M}abuchi energy
  asymptotics.
\newblock {\em Annals of Math.}, (175), 2012.

\bibitem{phillipon1991}
P.~Phillipon.
\newblock Sur des hauters alternatives {I}.
\newblock {\em Mathematische Annalen}, 289:255--283, 1991.

\bibitem{smirnov2004}
A.~V. Smirnov.
\newblock Classification of nearly closed orbits for the action of semisimple
  complex linear groups on the projective spaces.
\newblock In {\em Invariant theory in all characteristics}, volume~35 of {\em
  CRM Proc. Lecture Notes}, pages 251--257. Amer. Math. Soc., Providence, RI,
  2004.

\bibitem{smirnov2005}
A.~V. Smirnov.
\newblock Projective orbits of reductive groups, and {B}rion polytopes.
\newblock {\em Uspekhi Mat. Nauk}, 60(2(362)):147--148, 2005.

\bibitem{soule1991}
C.~Soul\'e.
\newblock G\'eometrie d'{A}rakelov et th\'eorie des nombres transcendants.
\newblock {\em Ast\'erisque}, 198-199-200:355--371, 1991.

\bibitem{sturmfels2017}
Bernd Sturmfels.
\newblock The {H}urwitz form of a projective variety.
\newblock {\em Journal of Symbolic Computation}, 79:186--196, 2017.

\bibitem{tianferm}
Gang Tian.
\newblock On {K}\"ahler-{E}instein metrics on certain {K}\"ahler manifolds with
  {$C\sb 1(M)>0$}.
\newblock {\em Invent. Math.}, 89(2):225--246, 1987.

\bibitem{tianberg}
Gang Tian.
\newblock On a set of polarized {K}\"ahler metrics on algebraic manifolds.
\newblock {\em J. Differential Geom.}, 32(1):99--130, 1990.

\bibitem{tian94}
Gang Tian.
\newblock The {$K$}-energy on hypersurfaces and stability.
\newblock {\em Comm. Anal. Geom.}, 2(2):239--265, 1994.

\bibitem{tian97}
Gang Tian.
\newblock K\"ahler-{E}instein metrics with positive scalar curvature.
\newblock {\em Invent. Math.}, 130(1):1--37, 1997.

\bibitem{tiansurv}
Gang Tian.
\newblock K\"ahler-{E}instein manifolds of positive scalar curvature.
\newblock In {\em Surveys in differential geometry: essays on Einstein
  manifolds}, Surv. Differ. Geom., VI, pages 67--82. Int. Press, Boston, MA,
  1999.

\bibitem{bottchrnfrms}
Gang Tian.
\newblock Bott-{C}hern forms and geometric stability.
\newblock {\em Discrete Contin. Dynam. Systems}, 6(1):211--220, 2000.

\bibitem{tianbook}
Gang Tian.
\newblock {\em Canonical metrics in {K}\"ahler geometry}.
\newblock Lectures in Mathematics ETH Z\"urich. Birkh\"auser Verlag, Basel,
  2000.
\newblock Notes taken by Meike Akveld.

\bibitem{tianCPAM}
Gang Tian.
\newblock K-stability and {K}\"ahler-{E}instien metrics.
\newblock {\em Comm. Pure and Appl. Math.}, 68:1085 -- 1156, 2015.

\bibitem{zhang}
Shouwu Zhang.
\newblock Heights and reductions of semi-stable varieties.
\newblock {\em Compositio Math.}, 104(1):77--105, 1996.

\end{thebibliography}
\end{document}